\documentclass[aos,preprint]{imsart}
\usepackage{mathrsfs}

\RequirePackage[OT1]{fontenc}
\usepackage[authoryear,round,sort&compress]{natbib}   
\usepackage{amsmath,amssymb,amsfonts}
\usepackage{amsbsy, amsthm,epsfig}
\usepackage{graphicx,rotating, tabularx, booktabs}
\usepackage{color}
\usepackage{pdflscape}
\usepackage{subfigure}
\usepackage{amsmath}
\usepackage{amsthm}
\usepackage{amssymb}
\usepackage{color}
\usepackage{lscape}
\usepackage{rotating}
\usepackage{caption2}
\usepackage{subfigure}
\usepackage{graphicx}
\usepackage{amssymb}
\usepackage{multirow}
\usepackage{epstopdf}

\startlocaldefs
\numberwithin{equation}{section}
\theoremstyle{plain}

\endlocaldefs
\def\be{\begin{equation}}
\def\ee{\end{equation}}
\def\bea{\begin{eqnarray}}
\def\eea{\end{eqnarray}}
\def\bd{\begin{displaymath}}
\def\ed{\end{displaymath}}
\def\bda{\begin{eqnarray*}}
\def\eda{\end{eqnarray*}}
\def\bsm{\begin{small}}
\def\esm{\end{small}}

\def\nn{\nonumber}

\def\ha1{\hat \beta_1}

\def\T{{ \mathrm{\scriptscriptstyle T} }}

\def\bsc{\begin{scriptsize}}
\def\esc{\end{scriptsize}}

\usepackage{amssymb}
\def\be{\begin{equation}}
\def\ee{\end{equation}}
\def\bea{\begin{eqnarray}}
\def\eea{\end{eqnarray}}
\def\bd{\begin{displaymath}}
\def\ed{\end{displaymath}}
\def\bda{\begin{eqnarray*}}
\def\eda{\end{eqnarray*}}

\def\nn{\nonumber}

\def\ha1{\hat \beta_1}

\def\bsc{\begin{scriptsize}}
\def\esc{\end{scriptsize}}

\usepackage{booktabs}
\usepackage{threeparttable}

\usepackage{accents}

\newcommand\munderbar[1]{%
  \underaccent{\bar}{#1}}

\def\sep{\mathrm{sep}}

\makeatletter
\newcommand{\figcaption}{\def\@captype{figure}\caption}
\newcommand{\tabcaption}{\def\@captype{table}\caption}
\makeatother

\def\s1{{ \mathrm{\scriptscriptstyle (1)} }}
\def\ss2{{ \mathrm{\scriptscriptstyle (2)} }}
\def\be{\begin{equation}}
\def\ee{\end{equation}}
\def\bea{\begin{eqnarray}}
\def\eea{\end{eqnarray}}

\newcommand{\E}{\rm E}

\newtheorem{tm}{Theorem}

\newtheorem{pn}{Proposition}

\theoremstyle{definition}

\newcommand{\bX}{{\mathbf X}}
\newcommand{\bY}{{\mathbf Y}}

\newcommand{\bSigma}{\boldsymbol{\Sigma}}

\newcommand{\bPsi} {\boldsymbol{\Psi}}

\newtheorem{innercustomgeneric}{\customgenericname}
\providecommand{\customgenericname}{}
\newcommand{\newcustomtheorem}[2]{%
  \newenvironment{#1}[1]
  {%
   \renewcommand\customgenericname{#2}%
   \renewcommand\theinnercustomgeneric{##1}%
   \innercustomgeneric
  }
  {\endinnercustomgeneric}
}

\newcustomtheorem{customthm}{Theorem}
\newcustomtheorem{customas}{Assumption}

\renewcommand{\theequation}{\thesection.\arabic{equation}}
\renewcommand{\thetable}{ \arabic{table}}

\begin{document}

\begin{frontmatter}
\title{Multi-level Thresholding Test for  High Dimensional 
Covariance Matrices}
\runtitle{Thresholding test for covariance}
\begin{aug}

\author{\fnms{Song Xi} \snm{Chen}$^{1}$
\ead[label=e1]{csx@gsm.pku.edu.cn}}
\address{Guanghua School of Management and\\
Center for Statistical Science\\
Peking University\\
Beijing, 100871, China\\
\printead{e1}\\
},
\author{\fnms{Bin} \snm{Guo}$^{2}$}
\ead[label=e2]{guobin@swufe.edu.cn}
\address{Center of Statistical Research and\\
School of Statistics\\
Southwestern University of Finance\\
~~and Economics\\
Chengdu, Sichuan, 611130, China\\
\printead{e2}} and \author{\fnms{Yumou} \snm{Qiu}$^{3}$
\ead[label=e3]{yumouqiu@iastate.edu}}
\address{Department of Statistics\\
Iowa State University\\
Ames, Iowa 50010, USA\\
\printead{e3}}

\runauthor{S. X. Chen, B. Guo and Y. Qiu}

\affiliation{Peking University$^1$,
Southwestern University of Finance and Economics$^2$ and Iowa State University$^3$}


\end{aug}

\begin{abstract}
We consider testing the equality of two high-dimensional covariance matrices
by carrying out a multi-level thresholding procedure, which is designed to detect sparse and faint differences between the covariances.
A novel $U$-statistic composition is developed to  establish the asymptotic distribution of the thresholding statistics in conjunction with the matrix blocking and the coupling techniques.
We propose a multi-thresholding test that is shown to be
powerful in detecting sparse and weak differences  between two covariance matrices. The test is shown to have attractive detection boundary and to  attain the optimal minimax rate in the signal strength under different regimes of high dimensionality and  the sparsity of the signal.
Simulation studies are conducted to demonstrate the utility of the proposed test.
\end{abstract}



\end{frontmatter}

\renewcommand{\theequation}{\thesection.\arabic{equation}}
\setcounter{equation}{0}
\section{Introduction}

Understanding the dependence among data components is an important goal in high-dimensional data analysis as different dependence structures lead to different inference procedures, for instance in the  Hotelling's test for the mean \citep{Hotelling_1931} and  Fisher's linear discriminant analysis, the pooled covariance estimate is used under the assumption of the same covariance matrix between the two samples.  For high dimensional data, the covariance matrices are utilized in the form of its inverse to enhance the signal strength in the innovated Higher Criticism test for high dimensional means \citep{Hall_Jin_2010} and in {Gaussian graphical models \citep{Liu_2013, Ren_2015}.
In genetic studies, covariances are widely used to understand the interactions among genes, to study functionally related genes \citep{Yi_2007}, and to construct and compare co-expression genetic networks \citep{Fuente_2010}. 

As a multivariate statistical procedure is likely constructed based on a specific dependence structure of the data,  testing for the equality of two covariance matrices ${\bf\Sigma}_1$ and ${\bf\Sigma}_2$ from two populations has been an enduring task.
\cite{John_1971}, \cite{Gupta_Giri_1973}, \cite{Nagao_1973} and \cite{Perlman_1980} presented studies under the conventional fixed dimensional setting; see \cite{Anderson_2003} for a comprehensive review.
The modern high-dimensional data have generated a renewed interest under the so-called ``large $p$, small $n$" paradigm.  For Gaussian data with the dimension $p$ and the sample size $n$ being the same order,
\cite{Schott_2007} and \cite{Srivastava_Yanagihara_2010} proposed two sample tests based on the distance measure $\|{\bf
\Sigma}_1-{\bf \Sigma}_2\|^2_F$, the squared Frobenius matrix norm between the two covariances.
\cite{Bai_2009} considered a corrected likelihood ratio test via the large dimensional random matrix theory.
For nonparametric settings without explicitly restricting $p$ and the sample sizes,
\cite{Li_Chen_2012} proposed an $\ell_{2}$-test based on a linear combination of $U$-statistics which is an unbiased estimator of $\|{\bf \Sigma}_1-{\bf \Sigma}_2\|^2_F$.
\cite{QC_2012} studied an $\ell_{2}$-test for the bandedness of a covariance.
\cite{Cai_Liu_Xia_2013} proposed a test based on the maximal standardized differences (an $\ell_{\max}$-type formulation) between the entries of two sample covariance matrices. 
\cite{Chang_2017} constructed a simulation based approach to approximate the distribution of the maximal statistics.
Studies have shown that
the $\ell_{2}$-tests are powerful for detecting dense and weak differences in the covariances, while the $\ell_{\max}$-formulation is powerful against sparse and strong signal settings.

Detecting rare and faint signals has attracted much attention in high-dimensional statistical inference.  The studies have been largely concentrated for the mean problems (\citealp{Fan_1996}; \citealp{Donoho_Jin_2004}; \citealp{Delaigle_2011_JRSSB}; \citealp{Zhong_Chen_Xu_2013}; \citealp{Qiu_Chen_Nettleton_2016}), while 
 studies 
 for covariance matrices are much less. 
\cite{AC_2012} investigated near optimal testing rules for detecting nonzero correlations in a one sample setting 
for Gaussian data with clustered nonzero signals. 

The aim of this paper is on enhancing the power performance in testing differences between two covariances when the differences are both sparse and faint, which is the most challenging setting for signal detection and brings about the issue of optimal detection boundary for covariance matrices. We introduce  thresholding  on the $\ell_2$-formulation of \cite{Li_Chen_2012} to remove those non-signal bearing entries of the covariances, which reduces the overall noise (variance) level of the test statistic and increases the signal to noise ratio for the testing problem.
The formulation may be viewed as a parallel development to  the  thresholding method for detecting differences in the means, for instance the  Higher Criticism (HC) test of
\cite{Donoho_Jin_2004, Hall_Jin_2010} and \cite{Delaigle_2011_JRSSB}, and the $\ell_2$-thresholding formulation in \cite{Fan_1996}, \cite{Zhong_Chen_Xu_2013} and \cite{Qiu_Chen_Nettleton_2016}.
However, comparing with the studies on the thresholding tests for the means, there is few work on the thresholding tests for covariance matrices beyond a discussion in \cite{Donoho_Jin_2015}, largely
 due to a difficulty {in} treating the dependence among the entries of the sample covariance matrices. 

To overcome the theoretical difficulty, 
we adopt  a matrix version of the blocking method to partition the matrix entries to big square blocks separated by small rectangular blocks.
The coupling technique is used to  construct an  equivalent U-statistic 
to the thresholding test statistic based on the covariance matrix block partition.
The equivalent U-statistic formulation allows establishing  the martingale central limit theorem \citep{Hall_Heyde_1980} for the asymptotic distribution of the test statistic.
A multi-thresholding test procedure is
 proposed to make the test adaptive to the unknown signal strength and sparsity.
Under the setting of rare and faint differences between the two covariances, the power of the proposed test is studied and its detection boundary is derived, which shows the benefits of the multi-thresholding over existing two sample covariance tests.

The paper is organized as follows. We introduce the setting of the 
covariance testing in Section 2. The thresholding statistic and the multi-level thresholding test are proposed in Sections 3 and 4, with its power  and detection boundary established in Section 5.
Simulation studies and discussions are presented in Sections
6 and 7, respectively.
Proofs and a real data analysis are relegated to the appendix and the supplementary material (SM).

\setcounter{equation}{0}
\section{Preliminary}

Suppose that there are two independent samples of  $p$-dimensional random vectors $\bX_{1},\dots,\bX_{n_1}
 \overset{i.i.d.}{\sim} \text{F}_1$ and $\bY_{1},\dots, \bY_{n_2} \overset{i.i.d.}{\sim} \text{F}_2$
drawn from two distributions $\text{F}_1$ and $\text{F}_2$, respectively,
where $\bX_{k}=(X_{k 1},\dots,X_{k p})^{\T}$,
$\bY_{k}=(Y_{k  1},\dots,Y_{k p})^{\T}$, $n_1$ and $n_2$ are the sample sizes, and ``i.i.d.'' stands for ``independent and identically distributed''.
Let ${\bf \mu}_1 = (\mu_{11}, \ldots, \mu_{1p})^{\T}$ and ${\bf\mu}_2 = (\mu_{21}, \ldots, \mu_{2p})^{\T}$ be the means of
$\text{F}_1$ and $\text{F}_2$, and
${\bf\Sigma}_1=(\sigma_{ij1})_{p\times p}$ and ${\bf
\Sigma}_2=(\sigma_{ij2})_{p\times p}$ be the covariance matrices of
$\text{F}_1$ and $\text{F}_2$,  respectively.
Let $\bPsi_{1} = (\rho_{ij1})_{p\times p}$ and $\bPsi_{2} = (\rho_{ij2})_{p\times p}$ be the corresponding correlation matrices.
We consider testing
\be H_0:{\bf \Sigma}_1={\bf \Sigma}_2  \quad \text{vs.} \quad
H_a:{\bf \Sigma}_1\neq{\bf \Sigma}_2 \label{H0} \ee
under a high-dimensional setting where $p \gg n_1,n_2$.

Let ${\bf \Delta} = {\bf \Sigma}_1 - {\bf \Sigma}_2 = (\delta_{ij})$ where $\delta_{ij} = \sigma_{ij1} - \sigma_{ij2}$ are component-wise differences between ${\bf \Sigma}_1$ and ${\bf \Sigma}_2$, 
 $q=p(p+1)/2$ be the number of distinct 
 parameters and
$n = n_1n_2 / (n_1 + n_2)$ be the effective sample size in the testing problem.

While hypothesis (\ref{H0}) offers all possible alternatives against the equality of the two covariances,
 we consider in this study a subset of the alternatives that constitutes the most challenging setting with the number of non-zero $\delta_{ij}$ %
being rare  and the magnitude of the non-zero $\delta_{ij}$ being faint; see
\cite{Donoho_Jin_2004} and \cite{Hall_Jin_2010}
for similar settings in the context for testing means.
Let $m_a$ denote the number of nonzero $\delta_{ij}$ for $i \leq j$. 
We assume a sparse setting such that  $m_a = \lfloor q^{(1-\beta)} \rfloor$ for a $\beta\in(1/2,1)$, where $\beta$ is the sparsity parameter and $\lfloor \cdot \rfloor$ is the integer truncation  function. 
We note that 
 $\beta \in (0,1/2]$ 
is the dense case under which the testing is easier. 

The faintness of signals is characterized by
\be
\delta_{ij} = \sqrt{2r_{0, ij}\log(q)/n} = \sqrt{4r_{0, ij}\log(p)/n}\{1 + o(1)\} ~~ \mbox{if~~ $\delta_{ij} \ne 0$}
\label{eq:Sstrength}\ee
for $r_{0, ij}>0$. 
As shown in Theorem \ref{tm3}, 
 $\sqrt{\log(p) / n}$ in (\ref{eq:Sstrength}) is the minimum rate for successful signal detection under the sparse setting. 
Specifically, our analysis focuses on  a special case of (\ref{H0}) such that
\be\begin{split}
&H_{0}: \delta_{ij} = 0 \mbox{ \ for all $1\leq i \leq j \leq p$ \ vs.} \\
&H_{a}: \mbox{ 
$m_a =\lfloor q^{(1-\beta)} \rfloor$ 
nonzero $\delta_{ij}$ with strength specified in (\ref{eq:Sstrength}).}
\end{split}
\label{SparseH1}\ee
Here, the signal strength 
 $r_{0, ij}$
together with $\beta \in (1/2,1)$ constitutes the rare and faint signal setting,
which has been used to evaluate tests on  
means and regression coefficients (\citealp{Donoho_Jin_2004}; \citealp{Hall_Jin_2010}; \citealp{Zhong_Chen_Xu_2013}; \citealp{Qiu_Chen_Nettleton_2016}). 
{Our proposed test 
 is designed to achieve high power under $H_a$ of (\ref{SparseH1}) that offers the most challenging setting for detecting unequal covariances, as shown in Theorem \ref{tm3}.}

Hypotheses (\ref{SparseH1}) are composite null versus composite alternative. Under the null, although  the two covariances
are the same, they  can take different values; and under the alternative no prior 
distribution is assumed on the location of the nonzero $\delta_{ij}$.
This is different from the simple null versus simple alternative setting of \cite{Donoho_Jin_2004}. 
The derivation of the optimal detection boundary for such composite hypotheses is more difficult as shown in the later analysis. 

Let $\{\pi_{\ell, p}\}_{\ell=1}^{p !}$ denote all possible permutations of $\{1, \dots, p\}$ 
and ${\bX}_{k}(\pi_{\ell, p})$ and ${\bY}_{k}(\pi_{\ell, p})$ be the reordering of ${\bX}_{k}$ and ${\bY}_{k}$ corresponding to a permutation $\pi_{\ell, p}$.
We assume that there is a permutation $\pi_{\ell_{*}, p}$ such that ${\bX}_{k}(\pi_{\ell_{*}, p})$ and ${\bY}_{k}(\pi_{\ell_{*}, p})$ are weakly dependent,  defined via the $\beta$-mixing { \citep{Bradley_2005}.
As the proposed statistic 
 in (\ref{SingleTest})} is of the $\ell_2$-type and is invariant to  permutations of 
 ${\bX}_{k}$ and ${\bY}_{k}$, there is no need to know  $\pi_{\ell_{*}, p}$. 

Let  $\{{\bX}_{k}\} = \{{\bX}_{k}(\pi_{\ell_{*}, p})\}$ and $\{{\bY}_{k}\} = \{{\bY}_{k}(\pi_{\ell_{*}, p})\}$ to simplify  notation. 
Let $\mathcal{F}_{m_a}^{m_b}({\bX}_{k})=\sigma\{X_{kj}: m_a \leq j \leq m_b\}$ and $\mathcal{F}_{m_a}^{m_b}({\bY}_{k})=\sigma\{Y_{kj}: m_a \leq j \leq m_b\}$ be the $\sigma$-fields generated by $\{{\bX}_{k}\}$ and $\{{\bY}_{k}\}$ for $1 \leq m_a \leq m_b \leq p$.
The $\beta$-mixing coefficients  
 are
$\zeta_{x, p}(h) = \sup_{1 \leq m \leq p - h} \zeta\{\mathcal{F}_{1}^{m}({\bX}_{k}), \mathcal{F}_{m+h}^{p}({\bX}_{k})\}$ and $\zeta_{y, p}(h) = \sup_{1 \leq m \leq p - h} \zeta\{\mathcal{F}_{1}^{m}({\bY}_{k}), \mathcal{F}_{m+h}^{p}({\bY}_{k})\}$ \citep{Bradley_2005},
where for two $\sigma$-fields $\mathcal{A}$ and $\mathcal{B}$, 
$$\zeta(\mathcal{A}, \mathcal{B}) = \frac{1}{2}\sup \sum_{l_1 = 1}^{u_1}\sum_{l_2 = 1}^{u_2}\big| P(A_{l_1} \cap B_{l_2}) - P(A_{l_1})P(B_{l_2}) \big|.$$
Here, the supremum is taken over all finite partitions $\{A_{l_1} \in \mathcal{A} \}_{l_1=1}^{u_1}$  and $\{B_{l_2} \in \mathcal{B}\}_{l_2=1}^{u_2}$ of the sample space, and 
 $u_1,u_2 \in \mathbb{Z}^{+}$, the set of positive integers.

Let $\bar{\bX}=\sum_{k=1}^{n_1}\bX_k/n_{1}$ and
$\bar{\bY}=\sum_{k=1}^{n_2}\bY_k/n_{2}$ be the two sample means
where
$\bar{\bX} = (\bar{X}_1, \dots, \bar{X}_p)^{\T}$ and $\bar{\bY} =
(\bar{Y}_1, \dots, \bar{Y}_p)^{\T}$.
Let
\bea
\widehat{{\bf\Sigma}}_1 &=& (\hat{\sigma}_{ij1}) \ = \  \frac{1}{n_1}\sum_{k=1}^{n_1}(\bX_{k}-\bar{\bX})(\bX_{k}-\bar{\bX})^{\T} \ \mbox{and } \nn \\
\widehat{{\bf\Sigma}}_2 &=& (\hat{\sigma}_{ij2}) \ = \ \frac{1}{n_2}\sum_{k=1}^{n_2}(\bY_{k}-\bar{\bY})(\bY_{k}-\bar{\bY})^{\T},  \nn
\eea
and $\kappa = \lim_{n_1,n_2 \to \infty} {n_1}/(n_1+n_2)$. 
Moreover, let $\theta_{ij1}=\text{var}\{(X_{ki}-\mu_{1i})(X_{kj}-\mu_{1j})\}$,
$\theta_{ij2}=\text{var}\{(Y_{ki}-\mu_{2i})(Y_{kj}-\mu_{2j})\}$;  
$\rho_{ij,lm}^{\s1} = \text{Cor}\{(X_{ki} - \mu_{1i})(X_{kj} - \mu_{1j}), (X_{kl} - \mu_{1l})(X_{km} - \mu_{1m})\}$, and $\rho_{ij,lm}^{\ss2} = \text{Cor}\{(Y_{ki} - \mu_{2i})(Y_{kj} - \mu_{2j}), (Y_{kl} - \mu_{2l})(Y_{km} - \mu_{2m})\}$.
Both $\theta_{ij1}$ and $\theta_{ij2}$ 
 can be estimated by
\bea
\hat{\theta}_{ij1} &=& \frac{1}{n_1}\sum_{k=1}^{n_1}\{(X_{ki}-\bar{X}_i)(X_{kj}-\bar{X}_j)-\hat{\sigma}_{ij1}\}^2 \  \mbox{and} \nn \\
\hat{\theta}_{ij2} &=& \frac{1}{n_2}\sum_{k=1}^{n_2}\{(Y_{ki}-\bar{Y}_i)(Y_{kj}-\bar{Y}_j)-\hat{\sigma}_{ij2}\}^2. \nn
\eea
%
As $\hat{\theta}_{ij1}/n_1+\hat{\theta}_{ij2}/n_2$ is  ratioly consistent to the variance of $\hat{\sigma}_{ij1}-\hat{\sigma}_{ij2}$,   we define a 
 standardized difference between $\hat{\sigma}_{ij1}$ and $\hat{\sigma}_{ij2}$ as
\[
M_{ij}=F_{ij}^{2} \mbox{ \ for \ } F_{ij}=\frac{\hat{\sigma}_{ij1}-\hat{\sigma}_{ij2}}{(\hat{\theta}_{ij1}/n_1+\hat{\theta}_{ij2}/n_2)^{1/2}}, \
 1\leq i\leq j \leq p.
\]

\cite{Cai_Liu_Xia_2013} proposed a maximum statistic $M_n=\max_{1\leq i\leq j\leq p} M_{ij}$ 
that targets at the largest signal between
${\bf\Sigma}_1$ and ${\bf\Sigma}_2$.
\cite{Li_Chen_2012} proposed an $\ell_2$-test that aims at $\|{\bf
\Sigma}_1-{\bf \Sigma}_2\|^2_F$. 
\cite{Donoho_Jin_2015} briefly 
 discussed the possibility of applying the Higher Criticism (HC) statistic for testing 
$H_0: {\bf\Sigma} = {\bf I}_p$ with Gaussian data. 

We are to 
propose a test by carrying out multi-level thresholding on $\{M_{ij}\}$ to filter out potential signals via an $\ell_2$-formulation, and 
show that such thresholding leads to
a more powerful test than both the maximum test and the $\ell_2$-type tests when the signals are rare and faint.

\setcounter{equation}{0}
\section{Thresholding statistics for covariance matrices}

By the moderate deviation result in Lemma 2 in
the SM, 
{under Assumptions \ref{as1} (or \ref{as1poly}), \ref{as2},  \ref{as3} and 
$H_{0}$ of (\ref{H0}),
 $P\big\{\max_{1\leq i\leq j\leq p} M_{ij}>4\log (p)\big\} \to 0$
 as $n, p \to \infty$. }
 This implies that a threshold level of $4\log (p)$ is asymptotically too large
under the null hypothesis,
and suggests a smaller threshold $\lambda_{p}(s)=4s\log (p)$ for a  thresholding parameter $s\in(0,1)$.
This leads to a thresholding statistic
\be\label{SingleTest}
T_n(s)=\sum_{1\leq i\leq j\leq p} M_{ij}\mathbb{I}\{M_{ij}>\lambda_{p}(s)\},
\ee
where $\mathbb{I}(\cdot)$ denotes the indicator function.

Statistic $T_n(s)$ removes those small standardized differences $M_{ij}$  between $\widehat{{\bf\Sigma}}_1$ and $\widehat{{\bf\Sigma}}_2$.
Compared with the $\ell_2$-statistic of \cite{Li_Chen_2012}, 
 $T_n(s)$ keeps only large $M_{ij}$ after filtering out the potentially insignificant ones.
By removing those smaller $M_{ij}$'s, the variance of $T_n(s)$ is much
reduced from 
that of \cite{Li_Chen_2012} 
which translates to a 
 larger power  
 as shown in the next section. 
Compared to the $\ell_{\max}$-test of \cite{Cai_Liu_Xia_2013} whose power is determined by
the maximum of $M_{ij}$, the thresholding
statistic not only uses the largest $M_{ij}$, but also all 
relatively large entries. 
This enhances the ability in detecting weak
signals as reflected in the power and the detection boundary in Section 5. 

Let $C$ be a positive constant whose value may change in the context.
For two real sequences $\{a_n\}$ and $\{b_n\}$, $a_n \sim b_n$ means that there are two positive constants $c_1$ and $c_2$ such that $c_1\leq a_n/b_n\leq c_2$ for all $n$.
We make the following assumptions in our analysis.
\begin{customas}{1A}
\label{as1} As $n\to \infty$, $p\to \infty$, $\log p \sim n^{\varpi}$ for a $\varpi \in (0, 1/5)$.
\end{customas}

\begin{customas}{1B}
\label{as1poly} As $n\to \infty$, $p\to \infty$,  $n \sim p^{\xi}$ for a $\xi \in (0, 2)$.
\end{customas}

\begin{customas}{2}
\label{as2} There exists a positive constant  $\tau$   such that
\bea
&\tau < \min_{1\leq i \leq p}\{\sigma_{ii1}, \sigma_{ii2}\} \leq \max_{1\leq i \leq p}\{\sigma_{ii1}, \sigma_{ii2}\} < \tau^{-1} \ \text{and} \label{Assum2-1} \\
&\min_{i,j}\{ \theta_{ij1}/(\sigma_{ii1}\sigma_{jj1}),\theta_{ij2}/(\sigma_{ii2}\sigma_{jj2})\} >\tau.\label{Assum2-2}
\eea
\end{customas}

\begin{customas}{3}
\label{as3} There exist positive constants $\eta$ and $C$ such that
for all $|t|<\eta$,
\[
E[\exp\{t(X_{ki}-\mu_{1i})^2\}]\leq
C~~\text{and}~~E[\exp\{t(Y_{ki}-\mu_{2i})^2\}]\leq C \quad
\text{for}~~i=1,\dots, p.
\]
\end{customas}

\begin{customas}{4}
\label{as4}
There exists a small positive constant $\rho_0$ such that 
\be
\max\{|\rho_{ij1}|, |\rho_{ij2}|\} < 1 - \rho_0 \mbox{ \ for any $i \neq j$}, 
\label{Assum4-1}\ee
and $\max\{|\rho_{ij,lm}^{\s1}|, |\rho_{ij,lm}^{\ss2}|\} < 1 - \rho_0$ for any $(i, j) \neq (l, m)$.
%
\end{customas}

\begin{customas}{5}
\label{as5} There is a permutation ($\pi_{\ell_{*}, p}$) of the 
 data sequences $\{X_{kj}\}_{j=1}^{p}$ and
$\{Y_{kj}\}_{j=1}^{p}$ 
 such that the permuted sequences  are $\beta$-mixing with the mixing coefficients satisfying $\max\{ \zeta_{x, p}(h), \zeta_{y, p}(h)\}\leq C
\gamma^{h}$ for a constant $\gamma \in (0, 1)$, any $p \in \mathbb{Z}^{+}$ and positive integer $h \leq p - 1$.
\end{customas}

Assumptions \ref{as1} and \ref{as1poly} specify 
 the exponential and polynomial growth rates of $p$ relative to $n$, respectively.
Assumption \ref{as2} prescribes that $\theta_{ij1}$ and $\theta_{ij2}$  are bounded away from zero to ensure the denominators of $M_{ij}$ being bounded away from zero with probability approaching 1. 
Assumption \ref{as3} assumes the distributions of 
 $X_{ki}$ and $Y_{ki}$ are sub-Gaussian. 
Sub-Gaussianity  is commonly assumed in high-dimensional literature \citep{BL_2008a, Xue_Ma_Zou_2012, Cai_Liu_Xia_2013}.
Assumption \ref{as4} regulates the correlations among variables in ${\bf X}_{k}$ and ${\bf Y}_{k}$, and subsequently the correlations among $\{ F_{ij}\}$  where $M_{ij} = F_{ij}^{2}$.

The $\beta$-mixing 
 Assumption \ref{as5} is made for the unknown variable permutation $\pi_{\ell_{*}, p}$. 
 Similar mixing conditions for the column-wise dependence were made in  \cite{Delaigle_2011_JRSSB} and \cite{Zhong_Chen_Xu_2013} for thresholding tests of means. 
If $\{X_{kj}\}_{j = 1}^{p}$ and $\{Y_{kj}\}_{j = 1}^{p}$ are both Markov chains (the vector sequence under the variable permutation), Theorem 3.3 in \cite{Bradley_2005} provides 
conditions for the processes being $\beta$-mixing.  If $\{X_{kj}\}_{j = 1}^{p}$ and $\{Y_{kj}\}_{j = 1}^{p}$ are linear processes with i.i.d.\ innovation processes $\{\epsilon_{x, kj}\}_{j = 1}^{p}$ and $\{\epsilon_{y, kj}\}_{j = 1}^{p}$, which include the ARMA processes as the special case, then they are $\beta$-mixing provided the innovation processes are absolutely continuous \citep{Mokkadem_1988}. The latter is  particularly weak.   
Under the Gaussian distribution,  any covariance that matches to the covariance of an ARMA process up to a permutation will be $\beta$-mixing.
Furthermore,  normally  distributed data with banded covariance or block diagonal covariance  
after certain variable permutation also satisfy this assumption.  The $\beta$-mixing coefficients are assumed to decay at an exponential rate in Assumption \ref{as5} to simplify proofs, while arithmetic rates can be entertained at the expense of more technical details.

There are implications of the $\beta$-mixing on $\sigma_{ij 1}$ and $\sigma_{ij 2}$ due to Davydov's inequality, which potentially restricts the signal level 
 $\delta_{ij} 
 = \sqrt{2 r_{0, ij} \log(q)/n}$.  
However, as the $\beta$-mixing is assumed for the unknown permutation $\pi_{\ell_{*}, p}$, which is likely not the ordering of the observed data, 
 the restriction  would be minimal.  In the unlikely event that the observed order of the data matches that under $\pi_{\ell_{*}, p}$, 
 the $\beta$-mixing  would imply that 
the signals would appear near the main diagonal of $\bSigma_1$ and $\bSigma_2$. 
However, as the power of the test is determined by the detectable  signal strength at or larger than the order $\sqrt{\log(q) / n}$, 
the effect of the $\beta$-mixing on the alternative hypothesis and the power 
 is limited as long as there exists a portion of differences 
with the standardized strength above the detection boundary $\rho^{*}(\beta,\xi)$ established in Propositions \ref{pn3} and \ref{pn4}.

Let $\mu_{{T}_n,0}(s)$ and $\sigma^2_{{T}_n,0}(s)$ be the mean and
variance of the thresholding statistic ${T}_n(s)$, respectively,  under 
 $H_0$. 
Let $\phi(\cdot)$ and $\bar{\Phi}(\cdot)$ be the density and survival functions of $N(0,1)$, 
 respectively.
Recall that $q=p(p+1)/2$.
The following proposition provides expansions of $\mu_{{T}_n,0}(s)$ and $\sigma^2_{{T}_n,0}(s)$.
\begin{pn}
\label{pn1} Under Assumptions \ref{as1} or \ref{as1poly} and Assumptions \ref{as2}-\ref{as5}, we have $\mu_{{T}_n,0}(s)=\tilde{\mu}_{{T}_n,0}(s)\{1+O(\lambda_{p}^{3/2}(s) n^{-1/2})\}$ where
\[
\tilde{\mu}_{{T}_n,0}(s) = q\{2\lambda_p^{1/2}(s)\phi(\lambda_p^{1/2}(s))+2\bar{\Phi}(\lambda_p^{1/2}(s))\}.
\]
In addition, under either \text{(i)} Assumption \ref{as1} with $s>1/2$ or
\text{(ii)} Assumption \ref{as1poly} with $s>1/2-\xi/4$,
$\sigma^2_{{T}_n,0}(s)=\tilde{\sigma}^2_{{T}_n,0}(s)\{1+o(1)\}$,
where
$
\tilde{\sigma}^2_{{T}_n,0}(s) = q[2\{\lambda_p^{3/2}(s)+3\lambda_p^{1/2}(s)\}\phi(\lambda_p^{1/2}(s))+6\bar{\Phi}(\lambda_p^{1/2}(s))].
$
\end{pn}

From Proposition \ref{pn1}, we see that the main orders $\tilde{\mu}_{{T}_n,0}(s)$ and $\tilde{\sigma}^2_{{T}_n,0}(s)$
of $\mu_{{T}_n,0}(s)$ and $\sigma^2_{{T}_n,0}(s)$ are known and are solely determined by $p$ and $s$, and hence
can be readily used to estimate the mean and variance of $T_{n}(s)$.
The smaller  order term $\lambda_{p}^{3/2}(s) n^{-1/2}$  in $\mu_{{T}_n,0}(s)$  is useful in analyzing the performance of the thresholding test as in (\ref{EstMeanDifference}) later.
Compared to the variance of the thresholding statistic on the means \citep{Zhong_Chen_Xu_2013}, the exact main order $\tilde{\sigma}^2_{{T}_n,0}(s)$ of $\sigma^2_{{T}_n,0}(s)$ requires a minimum bound on the threshold levels, which is due to the more complex dependence among $\{M_{ij}\mathbb{I}(M_{ij} > \lambda_p(s))\}$.
More discussion regarding this is provided after Theorem \ref{tm1}.

{Next, we derive the asymptotic distribution of  $T_{n}(s)$ at a given $s$.}
The testing for the covariances involves a more complex dependency structure than those in time series and spatial data. 
In particular, although the data vector is $\beta$-mixing under a permutation,
the vectorization of $(M_{ij})_{p \times p}$ is not necessarily a mixing sequence, 
as 
 the sample covariances in the same row or column are dependent since they share common segments of data.
As a result,  the conventional blocking plus the coupling approach \citep{Berbee_1979} for mixing series is insufficient to establish the asymptotic distribution of $T_{n}(s)$.

To tackle the challenge, we first use a combination of the matrix blocking,  as illustrated in Figure \ref{Fig_Demo} in Appendix, and the coupling method. 
Due to the circular dependence of the sample covariances, this only produces independence among the big matrix blocks 
 with none overlapping indices, and those matrix blocks that share common indices are still dependent.  
To respect this reality, we introduce a novel U-statistic representation (\ref{Ustat}),  
which allows the use of the martingale central limit theorem on the U-statistic representation to attain the asymptotic normality of $T_n(s)$.

\begin{tm}
\label{tm1} Suppose Assumptions \ref{as2}-\ref{as5} are satisfied. Then, under the $H_0$ of (\ref{H0}), and either \text{(i)} Assumption \ref{as1} with $s>1/2$ or
\text{(ii)} Assumption \ref{as1poly} with $s>1/2-\xi/4$, we have
\[
\sigma^{-1}_{{T}_n,0}(s)\{{T}_n(s)-\mu_{{T}_n,0}(s)\}\stackrel{d}\to
N(0,1) \quad \mbox{as~$n, p \to \infty$.}
\]
\end{tm}

As the dependence between $M_{i_1j_1}\mathbb{I}\{M_{i_1j_1}>\lambda_{p}(s)\}$ and $M_{i_2j_2}\mathbb{I}\{M_{i_2j_2}>\lambda_{p}(s)\}$ decreases as the threshold level $s$ increases, the restriction on $s$ in Theorem \ref{tm1} is to control the dependence among the thresholded sample covariances in $T_n(s)$.
Under Assumption \ref{as1poly} that prescribes the polynomial growth of $p$, 
the minimum threshold level that guarantees the Gaussian limit of $T_n(s)$ can be chosen as close to 0 as $\xi$ approaches 2.
Compared to the thresholding statistic on the means \citep{Zhong_Chen_Xu_2013}, 
the thresholding on the covariance matrices requires a larger threshold level in order to control the dependence among entries of the sample covariances. 

\section{Multi-Thresholding test}

{ To formulate the multi-thresholding test, we need to first construct a single level thresholding test based on Theorem \ref{tm1}.}
From Proposition \ref{pn1}, 
we note that $\tilde{\sigma}^2_{{T}_n,0}(s)/\sigma^2_{{T}_n,0}(s) \to 1$. 
Let $\hat{\mu}_{{T}_n,0}(s)$ be an estimate of ${\mu}_{{T}_n,0}(s)$ that satisfies
\be \label{sufficient} 
\hat{\mu}_{{T}_n,0}(s) -{\mu}_{{T}_n,0}(s) =o_p\{\tilde{\sigma}_{{T}_n,0}(s)\}.
\ee
By Slutsky's theorem, under (\ref{sufficient}), the conclusion of Theorem \ref{tm1} is still valid if $\mu_{{T}_n,0}(s)$ and $\sigma^2_{{T}_n,0}(s)$ are replaced by $\hat{\mu}_{{T}_n,0}(s)$ and $\tilde{\sigma}^2_{{T}_n,0}(s)$, respectively.
A natural choice of $\hat{\mu}_{{T}_n,0}(s)$ is the main order term $\tilde{\mu}_{{T}_n,0}(s)$ given in Proposition \ref{pn1}. According to the expansion of $\mu_{{T}_n,0}(s)$,
\be
\frac{\mu_{{T}_n,0}(s)-\tilde{\mu}_{{T}_n,0}(s)}{\tilde{\sigma}_{{T}_n,0}(s)}=O_p\{\lambda_p^{5/4}(s)p^{1-s}n^{-1/2}\},
\label{EstMeanDifference}\ee
which converges to zero under Assumption \ref{as1poly} and $s>1-\xi/2$.
Therefore, we reject the null hypothesis of (\ref{H0}) if
\be
{T}_n(s) > \tilde{\mu}_{{T}_n,0}(s) + z_{\alpha}{\tilde{\sigma}_{{T}_n,0}(s)},
\label{test} \ee where $z_{\alpha}$ is the upper $\alpha$ quantile
of $N(0, 1)$.
We would call (\ref{test}) the single level thresholding test, since it is based on a single $s$.

It is noted that  Condition (\ref{sufficient}) {is to  simplify} the  analysis on the thresholding statistic.
When estimators satisfying (\ref{sufficient})
are not available, 
 we may choose $\hat{\mu}_{T_n,0}(s)=\tilde{\mu}_{T_n,0}(s)$ while the lower threshold bound has to be chosen as $1 - \xi / 2$ to make (\ref{EstMeanDifference}) converge to 0.
More accurate estimator of $\mu_{T_{n}, 0}(s)$ can  be constructed
by establishing expansions for ${\mu}_{{T}_n,0}(s)$  and then correcting for the bias empirically. 
\cite{Delaigle_2011_JRSSB} found that more precise moderate deviation results can be derived for the bootstrap calibrated t-statistics, which provides more accurate estimator for the mean. 

Existing works (\citealp{Donoho_Jin_2004}; \citealp{Delaigle_2011_JRSSB})
have shown  that for detecting rare and faint signals in 
means,
a single level thresholding cannot make the testing procedure  adaptive to the unknown signal strength and sparsity.
However,  utilizing many  thresholding levels can capture the underlying sparse and faint signals. This is the path we take for the covariance testing problem.

Let $\mathcal{T}_n(s)=\tilde{\sigma}^{-1}_{{T}_n,0}(s)\{{T}_n(s)-\hat{\mu}_{{T}_n,0}(s)\}$ 
 be the standardization of ${T}_n(s)$.
We construct a multi-level thresholding statistic by  maximizing 
 $\mathcal{T}_n(s)$ over a range of thresholds. 
This is in the same spirit of the HC test of \cite{Donoho_Jin_2004} and the multi-thresholding test of \cite{Zhong_Chen_Xu_2013} for the means.
Define the multi-level thresholding statistic
\be\label{MultiThreshold}
{\mathcal{V}}_n(s_0)=\sup_{s\in \mathcal{S}(s_0)}\mathcal{T}_n(s),\ee
where $\mathcal{S}(s_0)=(s_0,1-\eta]$ for a lower bound $s_0$ and an arbitrarily small
positive constant $\eta$.
From Theorem \ref{tm1}, a choice of 
  $s_0$ is either $1/2$ or $1/2 - \xi/4$ depending on  $p$ having the exponential or  polynomial growth. 
Define
\be
\mathcal{S}_n(s_0)=\{s_{ij}: s_{ij}=M_{ij}/(4\log
(p))~\text{and}~s_0< s_{ij}\leq (1-\eta)\}.
\label{thresholdset}
\ee
Since both
$\hat{\mu}_{{T}_n,0}(s)$ and $\tilde{\sigma}_{{T}_n,0}(s)$ are
monotone decreasing, 
${\mathcal{V}}_n(s_0)$ can be attained on $\mathcal{S}_n(s_0)$ such that 
\be {\mathcal{V}}_n(s_0)=\sup_{s\in \mathcal{S}_n(s_0)}\mathcal{T}_n(s).
\label{sup_test} \ee
This reduces the computation 
to finite number of threshold levels. 
The asymptotic distribution of 
 ${\mathcal{V}}_n(s_0)$ 
is given in the following theorem. 

\begin{tm}
\label{tm2} Suppose conditions of Theorem \ref{tm1} and
(\ref{sufficient}) hold, under $H_0$ of (\ref{H0}),
\[
P\{ a(\log (p)){\mathcal{V}}_n(s_0)-b(\log (p), s_0, \eta)\leq x \} \to
\exp(-e^{-x}),
\]
where $a(y)=(2\log (y))^{1/2}$ and $b(y,s_0,\eta)=2\log
(y)+2^{-1}\log\log(y)-2^{-1}\log(\pi)\\+\log(1-s_0-\eta)$.
\end{tm}

This 
 leads to an asymptotic $\alpha$-level multi-thresholding test (MTT) that
rejects $H_0$ if \be {\mathcal{V}}_n(s_0)>\{ q_{\alpha}+b(\log
(p),s_0,\eta)\}/a(\log (p)),
\label{testproc}\ee
where $q_{\alpha}$ is the upper $\alpha$
quantile of the Gumbel distribution. 
The test is adaptive to the unknown signal strength  and sparsity as revealed in the next section.
However, the convergence of ${\mathcal{V}}_n(s_0)$ 
can be slow, which may cause certain degree of size distortion. To speed up the convergence, we will present a parametric bootstrap procedure with estimated covariances to approximate the null distribution of ${\mathcal{V}}_n(s_0)$ in Section 6.

\setcounter{equation}{0}
\section{Power and detection boundary}

We evaluate the power performance of the proposed thresholding test (\ref{testproc}) under the alternative hypothesis (\ref{SparseH1}) by deriving its detection boundary, and demonstrate its superiority over the $\ell_2$-type and $\ell_{\max}$-type tests.

A detection boundary is a phase transition diagram 
 in terms of the signal
strength  and sparsity parameters $(r, \beta)$.   
We first outline the notion in the context of testing for high dimensional means. 
\cite{Donoho_Jin_2004} considered testing hypotheses for means 
from $p$ independent $N(\mu_j, 1)$-distributed data  for 
\be
H_{0}^{(m)}: \mu_{j} = 0  \mbox{ for all $j$  \ vs. \ } H_{a}^{(m)}: \mu_1,\ldots,\mu_p \overset{i.i.d.}{\sim} (1 - \epsilon) \nu_{0} + \epsilon \nu_{\mu_{a}}
\label{eq:TestMean}\ee
where $\epsilon = p^{-\beta}$, $\mu_{a} = \sqrt{2 r \log (p)}$, $\beta \in (0, 1)$ and $r \in (0, 1)$,  $\nu_{0}$ and $\nu_{\mu_{a}}$ denote the point mass distributions at $0$ and $\mu_{a}$, respectively.  The high dimensionality is reflected by $p \to \infty$. 
Let
\be
\rho(\beta) = \left\{
  \begin{array}{l l}
   \max\{0, \beta - 1/2\} & \quad \text{if $0 < \beta \leq 3/4$,} \\
   (1 - \sqrt{1 - \beta})^{2} & \quad \text{if $3/4 < \beta < 1$.} \\
  \end{array} \right.
\label{eq: DBoundaryMean}\ee
\cite{Ingster_1997} showed that $r=\rho(\beta)$ is the optimal detection boundary for hypotheses (\ref{eq:TestMean}) under the Gaussian distributed data setting of \cite{Donoho_Jin_2004}, in the sense that
(i) for any test of hypothesis (\ref{eq:TestMean}), 
\be P(\mbox{Reject $H_{0}^{(m)} | H_{0}^{(m)}$}) + P(\mbox{Not reject $H_{0}^{(m)} | H_{a}^{(m)}$}) \to 1 \quad    
\mbox{ \ if \ $r < \rho(\beta)$;}
\label{eq: belowDB}\ee
 and (ii) there exists a  test  such that
\be P(\mbox{Reject $H_{0}^{(m)} | H_{0}^{(m)}$}) + P(\mbox{Not reject $H_{0}^{(m)} | H_{a}^{(m)}$}) \to 0 
\mbox{ \ if \ $r > \rho(\beta)$,}
\label{eq: aboveDB}\ee
as $n, p \to \infty$. 
\cite{Donoho_Jin_2004} showed that the HC test attains this detection boundary, and  thus is optimal. 
They also derived phase transition diagrams for non-Gaussian data.
See \cite{Zhong_Chen_Xu_2013} and \cite{Qiu_Chen_Nettleton_2016} in other constructions for testing means and regression coefficients that also have $r = \rho(\beta)$ as the detection boundary which is not necessarily optimal under nonparametric data distributions. 

Define the standardized signal strength 
\be
r_{ij} = r_{0, ij} / \{(1 - \kappa)\theta_{ij1} + \kappa \theta_{ij2}\} \quad \hbox{for} \quad \sigma_{ij1} \neq \sigma_{ij2},
\label{SSignal}\ee
by recognizing that the denominator 
 is the main order term of the variance of $\sqrt{n}(\hat{\sigma}_{ij1}-\hat{\sigma}_{ij2})$. 
Under Gaussian distributions, $\theta_{ij1} = \sigma_{ii1}\sigma_{jj1} + \sigma_{ij1}^2$ and $\theta_{ij2} = \sigma_{ii2}\sigma_{jj2} + \sigma_{ij2}^2$. 
Under the alternative hypothesis in (\ref{SparseH1}), since the difference between $\sigma_{ij1}$ and $\sigma_{ij2}$ is at most at the order $\sqrt{\log(p) / n}$, we have
$r_{ij} = r_{0, ij} / (\sigma_{ii1}\sigma_{jj1} + \sigma_{ij1}^{2}) \{1 + O(\sqrt{\log(p) / n})\}$. 

Define the maximal and minimal standardized signal strength
\be
\bar{r} = \max_{(i, j): \sigma_{ij1} \neq \sigma_{ij2}} r_{ij} \mbox{ \ and \ } 
\munderbar{r} = \min_{(i, j): \sigma_{ij1} \neq \sigma_{ij2}} r_{ij}. 
\label{MSSignal}\ee
%
Let 
\vspace{-0.5cm}
\be\begin{split}
\mathcal{C}(\beta, \bar{r}, \munderbar{r}) =&  \big\{({\bf \Sigma}_{1}, {\bf \Sigma}_{2}): 
\mbox{under $H_{a}$ of (\ref{SparseH1}) such that  $m_a= \lfloor q^{1-\beta} \rfloor$, maximal } \\
& \ \ \
\mbox{and minimal standardized signal strength are $\bar{r}$ and $\munderbar{r}$,} \\ 
& \ \ \
\mbox{respectively, and satisfy Assumptions \ref{as2}, \ref{as4} and \ref{as5}} \big\}
\end{split}\nn\ee
be the class of covariance matrices with sparse and weak differences. 
{For any $({\bf \Sigma}_{1}, {\bf \Sigma}_{2}) \in \mathcal{C}(\beta, \bar{r}, \munderbar{r})$, }
let $\mu_{{T}_n,1}(s)$ and $\sigma^2_{{T}_n,1}(s)$ be the mean and
variance of ${T}_n(s)$ under $H_a$ in (\ref{SparseH1}),
and let 
\be
\mbox{Power}_{n}({\bf \Sigma}_{1}, {\bf \Sigma}_{2}) = P\big[{\mathcal{V}}_n(s_0) > \{ q_{\alpha}+b(\log(p),s_0,\eta)\}/a(\log (p)) | {\bf \Sigma}_{1}, {\bf \Sigma}_{2}\big] \nn
\ee
be the power of the MTT in (\ref{testproc}). 
%
Put $\mbox{SNR}(s) = \frac{\mu_{T_n,1}(s) - \mu_{T_n,0}(s)}{{\sigma}_{{T}_n,1}(s)}$ be the signal to noise ratio under $H_{a}$ in (\ref{SparseH1}).
Note that 
\[{\mathcal{V}}_n(s_{0}) = \max_{s\in \mathcal{S}(s_0) } \frac{{\sigma}_{{T}_n,1}(s)}{\tilde{\sigma}_{{T}_n,0}(s)} \bigg\{\frac{{{T}}_n(s)-{\mu}_{{T}_n,1}(s)}{{\sigma}_{{T}_n,1}(s)} - \frac{\hat{\mu}_{{T}_n,0}(s) - {\mu}_{{T}_n,0}(s)}{{\sigma}_{{T}_n,1}(s)} +\text{SNR}(s) \bigg\}.
\]
Thus, the power of the MTT is critically determined by $\mbox{SNR}(s)$. 

The next proposition 
gives the mean and variance of $T_{n}(s)$ under $H_{a}$ of (\ref{SparseH1}) with the same standardized signal strength ${r}^{*}$, corresponding to the cases that $r_{ij} = \munderbar{r}$ for all $\sigma_{ij1} \neq \sigma_{ij2}$ (${r}^{*} = \munderbar{r}$)
and $r_{ij} = \bar{r}$ for all $\sigma_{ij1} \neq \sigma_{ij2}$ (${r}^{*} = \bar{r}$).
Let $L_p$ be a multi-$\log (p)$ term which may change in context.

\begin{pn}
\label{pn2} Under Assumptions \ref{as1} or \ref{as1poly}, \ref{as2}-\ref{as5} and $H_{a}$ in (\ref{SparseH1}) with $r_{ij} = {r}^{*}$ for all $\sigma_{ij1} \neq \sigma_{ij2}$, $\mu_{{T}_n,1}(s) = \mu_{T_n,0}(s)+\mu_{T_n,a}(s)$, where
\[
\mu_{{T}_n,a}(s) =  L_pq^{(1-\beta)}\mathbb{I}(s<{r}^{*})+L_pq^{(1-\beta)}p^{-2(\sqrt{s}-\sqrt{{r}^{*}})^2}\mathbb{I}(s>{r}^{*}).
\]
In addition, under either \text{(i)} Assumption \ref{as1} with $s>1/2$ or
\text{(ii)} Assumption \ref{as1poly} with $s>1/2-\xi/4$,  $\sigma^2_{{T}_n,1}(s)=L_pq^{(1-\beta)}p^{-2(\sqrt{s}-\sqrt{{r}^{*}})^2}\mathbb{I}(s>{r}^{*})+L_pq^{(1-\beta)}\mathbb{I}(s<{r}^{*})+L_pqp^{-2s}$.
\end{pn}

From Proposition \ref{pn2} via the maximal and minimal signal strength defined in (\ref{MSSignal}),
the detection boundary of the proposed MTT are established in Propositions \ref{pn3} and \ref{pn4} below.
As shown in the previous section, a lower threshold bound $s_{0}$ is needed to control the dependence among the entries of the sample covariance matrices. 
The restriction on the threshold levels leads to 
a slightly higher detection boundary as compared with that  given in (\ref{eq: DBoundaryMean}).  Before proceeding further,  let us define a family of detection boundaries indexed by  $\xi \in [0,2]$ that connects $p$ and $n$ via  $n \sim p^{\xi}$:
\be \label{detect_real} \rho^{*}(\beta,\xi)=\left\{\begin{matrix}
\frac{(\sqrt{4-2\xi}-\sqrt{6-8\beta-\xi})^2}{8},& 1/2<\beta\leq 5/8-\xi/16,\\
 \beta-1/2,& 5/8-\xi/16<\beta\leq 3/4,\\
(1-\sqrt{1-\beta})^2, & 3/4<\beta<1.
\end{matrix}\right.
\ee
It is noted that the phase diagrams  $\rho^{*}(\beta,\xi)$ are only defined over $\beta \in (1/2,1)$, the  sparse signal range. 
It  can be checked that $\rho^{*}(\beta,\xi) \ge \rho(\beta)$ for $\beta \in (1/2,1)$ and any $\xi \in [0,2]$.

The following  proposition considers the case of $n \sim p^{\xi}$ for $\xi \in (0, 2)$ as prescribed in Assumption \ref{as1poly},  a case considered in \cite{Delaigle_2011_JRSSB} in the context of mean  testing.

\begin{figure}[h]
\centering
\includegraphics[scale=0.46]{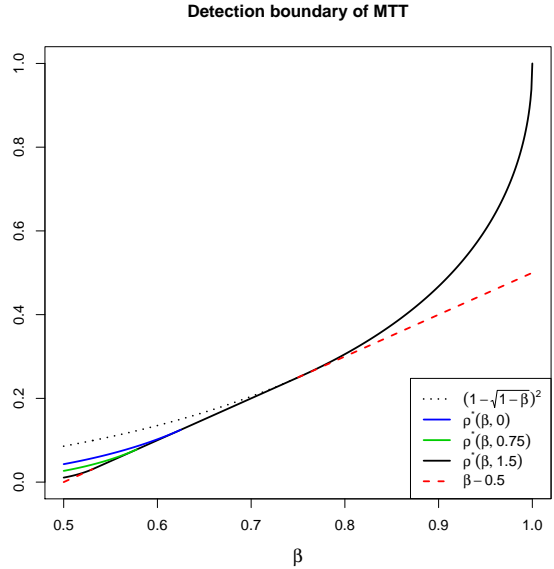}
\setlength{\abovecaptionskip}{0pt}
\setlength{\belowcaptionskip}{-5pt}
\caption{The detection boundary $\rho^{*}(\beta,\xi)$ in (\ref{detect_real}) of the proposed multi-level thresholding test 
 with $s_0 = 1/2 - \xi / 4$ for $\xi = 0, 0.75, 1.5$ and $n = p^{\xi}$, and the two pieces (in dashed and dotted curves) that constitute the optimal detection boundary $\rho(\beta)$ for testing means given in (\ref{eq: DBoundaryMean}). }
\label{Fig_DB}
\end{figure}

\begin{pn}
\label{pn3} 
Under Assumptions \ref{as1poly}, \ref{as2}-\ref{as5}, (\ref{sufficient}) and the alternative hypothesis (\ref{SparseH1}), for $s_0 = 1/2 - \xi / 4$, an arbitrarily small $\epsilon>0$, and a series of nominal sizes $\alpha_n=\bar{\Phi}((\log p)^{\epsilon})\to 0$, as $n, p \to \infty$,

\text{(i)} if $\munderbar{r} > \rho^{*}(\beta,\xi)$,
$\inf_{({\bf \Sigma}_{1}, {\bf \Sigma}_{2}) \in \mathcal{C}(\beta, \bar{r}, \munderbar{r})} \mbox{Power}_{n}({\bf \Sigma}_{1}, {\bf \Sigma}_{2}) \to 1$;

\text{(ii)} if $\bar{r} < \rho^{*}(\beta,\xi)$,
$\sup_{({\bf \Sigma}_{1}, {\bf \Sigma}_{2}) \in \mathcal{C}(\beta, \bar{r}, \munderbar{r})} \mbox{Power}_{n}({\bf \Sigma}_{1}, {\bf \Sigma}_{2}) \to 0$.
\end{pn}

{Proposition \ref{pn3} shows that the power of the proposed MTT over the class $\mathcal{C}(\beta, \bar{r}, \munderbar{r})$ is determined by $\beta$,  and the minimum and maximum standardized signal strength. More importantly, $\rho^{*}(\beta,\xi)$  in (\ref{detect_real}) is the detection boundary of the MTT.}
The power converges to 1 if $\munderbar{r}$ is above this boundary, and diminishes to 0 if $\bar{r}$ is below it. 
The detection boundaries $\rho^{*}(\beta,\xi)$ are displayed in Figure \ref{Fig_DB} for three values of $\xi$. 
Note that $\rho(\beta)$ in (\ref{eq: DBoundaryMean}) is the detection boundary of the MTT for $s_0 = 0$ that corresponds to $\xi=2$, which is the lowest one in the family.
It can be shown that $\rho^{*}(\beta,\xi)$ approaches to $\rho(\beta)$ as $\xi \to 2$; namely if $n \sim p^2$, we have $\rho^{*}(\beta,2) = \rho(\beta)$, which is the optimal detection boundary for testing the means with uncorrelated Gaussian data.
Restricting $s \geq s_0= 1/2 - \xi / 4$ 
 elevates the detection boundary $\rho^{*}(\beta,\xi)$ of the proposed MTT for $1/2<\beta\leq 5/8-\xi/16$ as a price for controlling the size of the test.
Similar results on the influence of the lower threshold bound on testing means were given in \cite{Delaigle_2011_JRSSB}. 

The following proposition shows that 
$\rho^{*}(\beta,0)$  is the   detection boundary when dimension $p$ grows exponentially fast with $n$, which can be viewed as a degenerated polynomial growth case with $\xi = 0$.

\begin{pn}
\label{pn4} Under Assumption \ref{as1}, \ref{as2}-\ref{as5}, (\ref{sufficient}) and the alternative hypothesis (\ref{SparseH1}), for $s_0 = 1/2$, an arbitrarily small $\epsilon>0$, and a series of nominal sizes $\alpha_n=\bar{\Phi}((\log p)^{\epsilon})\to 0$, as $n, p \to \infty$,

\text{(i)} if $\munderbar{r} > \rho^{*}(\beta,0)$,
$\inf_{({\bf \Sigma}_{1}, {\bf \Sigma}_{2}) \in \mathcal{C}(\beta, \bar{r}, \munderbar{r})} \mbox{Power}_{n}({\bf \Sigma}_{1}, {\bf \Sigma}_{2}) \to 1$;

\text{(ii)} if $\bar{r} < \rho^{*}(\beta,0)$,
$\sup_{({\bf \Sigma}_{1}, {\bf \Sigma}_{2}) \in \mathcal{C}(\beta, \bar{r}, \munderbar{r})} \mbox{Power}_{n}({\bf \Sigma}_{1}, {\bf \Sigma}_{2}) \to 0$.
\end{pn}

As $\rho^{*}(\beta,0) \ge \rho^{*}(\beta,\xi)$ for any $\xi \in (0,2]$, the result also shows that a higher growth rate of $p$ leads to a higher detection boundary that may be viewed as a sacrifice of the power due to the higher dimensionality.

From \cite{Cai_Liu_Xia_2013}, the power of the $\ell_{\max}$-test converges to 1 if
$$\max_{1 \leq i \leq j \leq p} \frac{|\sigma_{ij1} - \sigma_{ij2}|}{(\theta_{ij1} / n_1 + \theta_{ij2} / n_2)^{1/2}} > 4 \sqrt{\log p},$$
which is equivalent to $\bar{r} > 4$ in our context. Hence, the signal strength required by the $\ell_{\max}$-test is stronger than that  $r_{ij} \in (0, 1)$ required by
the MTT in this paper. 
Also, 
 the $\ell_2$-test of \cite{Li_Chen_2012} does not have non-trivial power
 for $\beta > 1/2$.
Hence, the proposed MTT is more powerful than both the $\ell_2$-tests and $\ell_{\max}$-tests in detecting sparse and weak signals. 

Propositions \ref{pn3} and \ref{pn4} indicate that 
the MTT can detect  the differences between the unequal covariances in ${\bf \Sigma}_{1}$ and ${\bf \Sigma}_{2}$ at the order of $c_{a}\sqrt{\log(p) / n}$ for some positive constant $c_{a}$. 
We are to show that the order  $\sqrt{\log(p) / n}$ is minimax optimal.

Let $\mathcal{W}_{\alpha}$ be the collection of all $\alpha$-level tests for  hypotheses (\ref{H0})
under Gaussian distributions and Assumptions \ref{as2}, \ref{as4} and \ref{as5},
namely, $P(W_{\alpha} = 1 | H_{0}) \leq \alpha$
for any $W_{\alpha} \in \mathcal{W}_{\alpha}$. 
Note that (\ref{Assum2-1}) and (\ref{Assum4-1}) are sufficient conditions for (\ref{Assum2-2}) and the second part of Assumption \ref{as4} under Gaussian Distributions, respectively.
Define a class of covariance matrices with the differences being at least of order $\{\log(p) / n\}^{1/2}$:
\bea
\underline{\mathcal{C}}(\beta, c) &=& \big\{({\bf \Sigma}_{1}, {\bf \Sigma}_{2}): 
\mbox{under $H_{a}$ of (\ref{SparseH1}) such that}
\ m_{a} = \lfloor q^{1-\beta} \rfloor, r_{0, ij} \geq c \nn \\ 
&& \ \ \ 
\mbox{for all $\sigma_{ij1} \neq \sigma_{ij2}$, and satisfy Assumptions \ref{as2}, \ref{as4} and \ref{as5}} \big\}. \nn
\eea
Having Assumptions \ref{as4} and \ref{as5} in $\underline{\mathcal{C}}(\beta, c_{0})$ and $\mathcal{W}_{\alpha}$ is for comparing the power performance of the MTT with the minimax rate.
Comparing with the covariance class $\mathcal{C}(\beta, \bar{r}, \munderbar{r})$, $\underline{\mathcal{C}}(\beta, c)$ has no constraint on the maximal signal strength.
For Gaussian data, $\theta_{ij1}, \theta_{ij2} \leq 2 \tau^{-2}$ where $\tau$ specifies the bounds in (\ref{Assum2-1}).
Thus,  the standardized signal strength $r_{ij} \geq c\tau^{2} / 2$ for all 
$\sigma_{ij1} \neq \sigma_{ij2}$.
For the MTT, from Propositions \ref{pn3} and \ref{pn4},
$\inf_{({\bf \Sigma}_{1}, {\bf \Sigma}_{2}) \in \underline{\mathcal{C}}(\beta, c)} \mbox{Power}_{n}({\bf \Sigma}_{1}, {\bf \Sigma}_{2}) \to 1$ as $n, p \to \infty$ for a large constant $c$.
The following theorem shows that the lower bound $\{\log(p) / n\}^{1/2}$ for the signal in $\underline{\mathcal{C}}(\beta, c)$ is the optimal rate, namely there is no $\alpha$-level test that can distinguish $H_{a}$ from $H_{0}$ in (\ref{SparseH1}) with probability approaching 1 uniformly over the class $\underline{\mathcal{C}}(\beta, c_{0})$ for some $c_{0} > 0$.

\begin{tm}
\label{tm3} For the Gaussian distributed data, under Assumptions \ref{as1poly}, \ref{as2}, \ref{as4} and \ref{as5}, for any $\tau > 0$, $0 < \omega < 1 - \alpha$ and $\max\{2/3, (3 - \xi) / 4\} < \beta < 1$, there exists a constant $c_{0} > 0$ such that,  as $n, p \to \infty$,
\[
\sup_{W_{\alpha} \in \mathcal{W}_{\alpha}} \inf_{({\bf \Sigma}_{1}, {\bf \Sigma}_{2}) \in \underline{\mathcal{C}}(\beta, c_{0})} P(W_{\alpha} = 1) \leq 1 - \omega.
\]
\end{tm}

As Propositions 3 and 4 have shown  that the proposed MTT can detect signals at the rate of $\{\log(p) / n\}^{1/2}$ for $\beta > 1/2$, the MTT test is at least minimax rate optimal for $\beta > \max\{2/3, (3 - \xi) / 4\}$.
Compared to Theorem 4 of \cite{Cai_Liu_Xia_2013}, by studying the alternative structures in  $\underline{\mathcal{C}}(\beta, c_{0})$,  we extend the minimax result from the highly sparse signal regime $3/4 < \beta < 1$ to $\max\{2/3, (3 - \xi) / 4\} < \beta < 1$,
 which offers a wider range of the signal sparsity.
The optimality under $1/2 < \beta \leq \max\{2/3, (3 - \xi) / 4\}$ requires  investigation in a separate effort.


Obtaining the lower and upper bounds of the detectable signal strength at the rate $\sqrt{\log(p) / n}$ requires  more sophisticated derivation.
These two bounds could be the same under certain conditions for testing one-sample covariances.
However, for the two-sample test, the lower and upper bounds may  not match.
This is due to the composite null hypothesis in (\ref{SparseH1}). 
{More discussion on this issue is given in Section 7.}

\setcounter{equation}{0}

\section{Simulation Results}

We report results from simulation experiments which were designed to evaluate the performances of the
proposed two-sample MTT under high dimensionality with sparse and faint signals.
We also compare the proposed test with the tests
in 
\cite{Srivastava_Yanagihara_2010} (SY),
\cite{Li_Chen_2012} (LC) and \cite{Cai_Liu_Xia_2013} (CLX).

In the simulation studies, the two random samples
$\{\bX_k\}_{k=1}^{n_1}$ and $\{\bY_k\}_{k=1}^{n_2}$ were respectively
generated from
\be
\bX_k={\bf \Sigma}_1^{\frac{1}{2}}{\bf Z}_{1k} \quad \text{and}\quad
\bY_k={\bf \Sigma}_2^{\frac{1}{2}}{\bf Z}_{2k}, \label{DataGenerate}
\ee
where $\{{\bf Z}_{1k}\}$ and $\{{\bf Z}_{2k}\}$ are i.i.d.\ random vectors from a common population.
We considered two distributions for the innovation vectors ${\bf Z}_{1k}$ and ${\bf Z}_{2k}$:
(i) $N(0,{\bf I}_p)$; (ii) Gamma distribution where  components of ${\bf Z}_{1k}$ and ${\bf Z}_{2k}$
were i.i.d.\ standardized Gamma(4,2) with mean 0 and variance 1.
To design the covariances ${\bf \Sigma}_1$ and ${\bf \Sigma}_2$, let ${\bf \Sigma}_1^{(0)}={\bf{D}}_0^{\frac{1}{2}}{\bf
\Sigma}^{(*)}{\bf{D}}_0^{\frac{1}{2}}$, where ${\bf{D}}_0 = \mbox{diag}(d_{1}, \ldots, d_{p})$ 
 with elements generated according to the uniform distribution $\text{U}(0.1,1)$, and ${\bf
\Sigma}^{(*)}=(\sigma_{ij}^{*})$ was a positive definite correlation matrix.
Once generated, ${\bf{D}}_0$  was held fixed throughout the simulation. 
The following two designs of ${\bf \Sigma}^{(*)}$ were considered in the simulation:
%
%
\bea
&\text{Design 1:}& \sigma_{ij}^{*} = 0.4^{|i-j|}; \label{Model1}
\\
&\text{Design 2:}& \sigma_{ij}^{*} = 0.5 \mathbb{I}(i = j) + 0.5 \mathbb{I}(i, j \in [4 k_0 - 3, 4 k_0] ). \label{Model2}
\eea
for $k_0 = 1, \ldots, \lfloor p/4\rfloor$.
Design 1 has an auto-regressive structure and Design 2 is block diagonal with block size 4.
Matrix ${\bf{D}}_0$ created heterogeneity for different dimensions of the data.

To generate scenarios of sparse and weak signals under the alternative hypothesis, we chose
\begin{equation}
{\bf \Sigma}_1^{(\star)} = {\bf \Sigma}_1^{(0)} + \epsilon_c {\bf I_p} \quad  \hbox{and} \quad
{\bf \Sigma}_2^{(\star)} = {\bf \Sigma}_1^{(0)} + {\bf U} + \epsilon_c {\bf I_p},
\label{Sigma_simu}
\end{equation}
{where ${\bf U} = (u_{kl})_{p\times p}$ is a banded symmetric matrix
and $\epsilon_c$ is a positive number to guarantee the positive definiteness of ${\bf \Sigma}_2^{(\star)}$.
Specifically, let $k_0 = \lfloor m_p / p \rfloor$, where $m_p = \lfloor q^{1-\beta} / 2\rfloor$ is the number of distinct pairs with nonzero $u_{kl}$.
Let $u_{l + k_0 + 1\, l} = u_{l\, l + k_0 + 1} = \sqrt{4r\log p/n}$ for $l = 1, \ldots, k_1$ and $k_1 = m_p - p k_0 + k_0(k_0 + 1) / 2$, and let $u_{kl} = \sqrt{4r\log p/n}$ for $|k - l| \leq k_0$ and $k \neq l$ if $k_0 \geq 1$.}
Set $\epsilon_c = | \min\{\lambda_{\min}({\bf\Sigma}_{1}^{(0)} + {\bf U}), 0\} | +0.05$, where $\lambda_{\min}(A)$ denotes the minimum eigenvalue of a matrix $A$.
Since $\epsilon_c > 0$  
and $\lambda_{\min}({\bf
\Sigma}_2^{(\star)}) \geq \lambda_{\min}({\bf\Sigma}_{1}^{(0)} + {\bf U}) + \epsilon_c > 0$, both ${\bf \Sigma}_1^{(\star)}$ and ${\bf \Sigma}_2^{(\star)}$ were positive definite under both Designs 1 and 2.
{Under the null hypothesis, we chose ${\bf \Sigma}_1 = {\bf \Sigma}_2 = {\bf \Sigma}_1^{(0)}$ in (\ref{DataGenerate}),
while under the alternative hypothesis, ${\bf \Sigma}_1 = {\bf \Sigma}_1^{(\star)}$ and ${\bf \Sigma}_2 = {\bf \Sigma}_2^{(\star)}$.

The simulated data were generated as a reordering of $\bX_k$ and $\bY_k$ from (\ref{DataGenerate})
according to a randomly selected permutation $\pi_p$ of $\{1, \ldots, p\}$. Once $\pi_p$ was generated, it was held fixed throughout the simulation.}
To mimic the regime of sparse and faint signals, we generated a set of $\beta$ and $r$ values.
First, we fixed $\beta=0.6$ and set $r=0.1, 0.2, \ldots, 1$ to create different signal strengths utilized in the simulation results shown in Figure \ref{Fig_DGP_1_r}. Then, $r=0.6$ was fixed while $\beta$ was varied from $0.3$ to $0.9$ to show the impacts of sparsity levels on the tests in Figure \ref{Fig_DGP_2_r}.
We chose the
sample sizes $(n_1,n_2)$  as $(60,60)$, $(80,80)$, $(100,100)$ and  $(120,120)$ respectively, and the corresponding dimensions $p=175, 277, 396$ and  $530$ according to $p=\lfloor 0.25n_1^{1.6}\rfloor$.
We set $s_0 =0.5$ according to Theorem \ref{tm1} and the discussion following (\ref{MultiThreshold}),
and $\eta$  was chosen as 0.05 in (\ref{thresholdset}).
{We chose $\hat{\mu}_{{T}_n,0}(s) = \tilde{\mu}_{{T}_n,0}(s)$.}
The process was replicated 500 times for each setting of the simulation.

Since the convergence of $\mathcal{V}_n(s_0)$ to the Gumbel distribution given in (\ref{testproc}) can be slow when the sample size was small, we employed a bootstrap procedure in conjunction with a consistent covariance estimator proposed by \cite{Rothman_2012}, 
which ensures the positive definiteness of the estimated covariance.
Since ${\bf\Sigma}_1={\bf\Sigma}_2$ under the null hypothesis,
the two samples $\{\bX_k\}_{k=1}^{n_1}$ and $\{\bY_k\}_{k=1}^{n_2}$ were pooled together to
estimate ${\bf\Sigma}_1$. Denote the estimator of \cite{Rothman_2012} as $\widehat{{\bf \Sigma}}$.
For the $b$-th bootstrap resample,
we drew
$n_1$ samples of
${\bX}^{*}$ and $n_2$ samples of ${\bY}^{*}$
independently  from $N(0,\widehat{{\bf \Sigma}})$.
Then, 
the bootstrap test statistic $\mathcal{V}_n^{*(b)}(s_0)$ was obtained based on ${\bX}^{*}$ and  ${\bY}^{*}$.
This procedure was repeated $B=250$ times to obtain the bootstrap sample of the proposed multi-thresholding statistic
$\{\mathcal{V}_n^{*(1)}(s_0),\ldots,\mathcal{V}_n^{*(B)}(s_0)\}$ under the null hypothesis.
The bootstrap empirical null
distribution of the proposed statistic was $\widehat{F}_0(x) = \frac{1}{B}\sum_{b = 1}^{B} \mathbb{I}\{\mathcal{V}_n^{*(b)}(s_0) \leq x\}$ and the bootstrap p-value was $1 - \widehat{F}_0(\mathcal{V}_n(s_0))$, where $\mathcal{V}_n(s_0)$ was the multi-thresholding statistic from the original sample.
We reject the null hypothesis if this p-value is smaller than the nominal significant level $\alpha=0.05$.
The validity of the bootstrap approximation can be justified in two key steps. First of all, if we generate the ``parametric bootstrap samples'' from the two normal distributions with the true population covariance matrices,
by Theorems 1 and 2, the bootstrap version of the single thresholding and multi-thresholding statistics will have the same limiting Gaussian  distribution and the extreme value distribution, respectively.
Secondly, we can replace the true covariance above by a consistently estimated  covariance matrix $\hat{{\bf \Sigma}}$ \citep{Rothman_2012}, which is {positive definite}.
The justification of the bootstrap procedure can be made by showing the consistency of $\hat{{\bf \Sigma}}$ by extending the results in \cite{Rothman_2012}.

Table\ref{size_gaussian} 
reports the empirical sizes of the proposed multi-thresholding test using the limiting Gumbel distribution for the critical value 
 (denoted as MTT) and the bootstrap calibration procedure described above (MTT-BT), together with three existing methods, with the nominal level 0.05, and the Gaussian and Gamma distributed random vectors, respectively. 
We observe that the MTT based on 
 the asymptotic distribution exhibited some size distortion when the sample size was small. 
However, with the increase of the sample size, the sizes of MTT became closer to the nominal level.
At the meantime,  the CLX and SY tests also experienced some size distortion under the Gamma scenario in smaller samples.
It is observed that the proposed multi-thresholding test with the bootstrap calibration (MTT-BT)  performed consistently well under all the scenarios with accurate empirical sizes.
This shows that the bootstrap distribution 
 offered more accurate approximation than the limiting Gumbel distribution to the distribution of the test statistic $\mathcal{V}_n(s_0)$ under the null hypothesis.

\newcolumntype{C}[1]{>{\centering\arraybackslash}p{#1}}

\begin{table}[!htb]
\caption{Empirical sizes for the tests  of
\cite{Srivastava_Yanagihara_2010} (SY),
\cite{Li_Chen_2012} (LC), \cite{Cai_Liu_Xia_2013} (CLX) and the proposed multi-level thresholding test based on the limiting distribution calibration in (\ref{testproc}) (MTT) and the bootstrap calibration (MTT-BT) for Designs 1 and 2 under the Gaussian and Gamma distributions with the nominal level of $5\%$. } \label{size_gaussian} 
\begin{center}
\begin{tabular}{cc|C{1.4cm}C{1.4cm}C{1.4cm}C{1.4cm}C{1.4cm}}
\hline \hline
  $ p $      & $(n_1, n_2)$ & SY & LC & CLX & MTT & MTT-BT \\
   \hline
             &              &\multicolumn{5}{c}{Gaussian Design 1} \\
\hline
   175    &  (60, 60)      &   0.048    &   0.058  &   0.054    & 0.088 & 0.058     \\
   \cline{3-7}
   277    &  (80, 80)      &  0.052    &   0.052  &    0.058 & 0.064 &  0.056
     \\
   \cline{3-7}
   396    &  (100, 100)      &   0.042 & 0.046 & 0.058 & 0.064 &  0.054  \\
   \cline{3-7}
   530    &  (120, 120)      &    0.056    &   0.048 &    0.050   & 0.056 & 0.046   \\
\hline
             &              &\multicolumn{5}{c}{Gaussian Design 2} \\
   \hline
   175    &  (60, 60)      &   0.060 & 0.048 & 0.052 & 0.094 & 0.048    \\
   \cline{3-7}
   277    &  (80, 80)      &   0.040 & 0.060 & 0.040 & 0.064 & 0.052    \\
   \cline{3-7}
   396   &  (100, 100)      &  0.052 & 0.042 & 0.044 & 0.090 &   0.048   \\
   \cline{3-7}
   530    &  (120, 120)      &  0.050 &0.046&0.044& 0.060 & 0.054  \\

\hline
             &              &\multicolumn{5}{c}{Gamma Design 1} \\
   \hline
   175    &  (60, 60)      & 0.046 & 0.060 & 0.066 & 0.110  & 0.056   \\
   \cline{3-7}
   277  &  (80, 80)      &  0.060 &	0.050	& 0.044  & 0.076 & 0.044   \\
   \cline{3-7}
   396    &  (100, 100)     &   0.046 &	0.052 &	0.046&    0.066   & 0.054  \\
   \cline{3-7}
   530    &  (120, 120)      &   0.060 & 0.056  & 0.048 & 0.060 &  0.048 \\

\hline
             &              &\multicolumn{5}{c}{Gamma Design 2} \\
   \hline
   175    &  (60, 60)      &  0.070 & 0.056 & 0.066 & 0.108    & 0.056 \\
   \cline{3-7}
   277   &  (80, 80)      &   0.060 &	0.058 &	0.068  & 0.112 & 0.044  \\
   \cline{3-7}
   396    &  (100, 100)     &  0.060 & 0.050 & 0.044 & 0.068 &  0.046  \\
   \cline{3-7}
   530    &  (120, 120)      & 0.054 & 0.056 & 0.048& 0.056 &  0.048 \\
\hline
\hline

\end{tabular}
\end{center}
\end{table}

Figure \ref{Fig_DGP_1_r} displays the
empirical powers with respect to different signal strengths $r$ for covariance matrix
Designs 1 and 2 with $n_1=n_2=80$, $p=277$ and $n_1=n_2=100$, $p=396$ under the Gaussian  distribution, respectively. Figure \ref{Fig_DGP_2_r} reports the empirical powers under different sparsity  ($\beta$) levels when the signal strength $r$ was fixed at 0.6.  Simulation results on the powers under the Gamma distribution are available in the SM.
It is noted that at $\beta=0.6$,  there were only 68 and 90 unequal entries  between {the upper triangles of} ${\bf \Sigma}_1$ and ${\bf \Sigma}_2$ among a total of $q=38503$ and $78606$ unique entries for $p=277$ and $396$, respectively.
To make the powers comparable for different methods, we adjusted the critical values of the tests by their respective empirical null distributions so that the actual sizes were approximately equal to the nominal level 5\%.
Due to the size adjustment, the MTT based on the limiting distribution and the MTT-BT based on the bootstrap calibration had the same test statistic, and hence the same power. Here, we only reported the numerical power results for the MTT-BT.

Figure  \ref{Fig_DGP_1_r} reveals  that the power of the proposed MTT-BT was the highest among all the tests under all the scenarios. Even though the powers  of other tests improved as the signal strength $r$ was increased, the proposed MTT-BT maintained a lead over the whole range of $r \in [0.1, 1]$.  The extra power advantage of the MTT-BT over the other three tests got larger as the signal strength $r$ increased.
We observe from Figure \ref{Fig_DGP_2_r} that the proposed test also had the  highest empirical power across the range of $\beta$.  The powers of the MTT-BT at the high sparsity level ($\beta \ge 0.7$) were higher than those of the CLX test. The latter test is known for doing well in the power when the signal was  sparse. We take this as an empirical confirmation to the attractive detection boundary of the proposed {MTT} established in the theoretical analysis reported in Section 5.
The monotone decrease pattern in the power profile of the four tests reflected the reality of reduction in the number of signals as $\beta$ was increased.
It is noted that the two $\ell_2$ norm based tests SY and LC are known to have good powers when the signals are dense, i.e. $\beta\leq 0.5$.
This was well reflected  in Figure \ref{Fig_DGP_2_r}  indicating the two tests had comparable powers to the MTT-BT when $\beta =0.3$ and $0.4$.
However, after $\beta$ was larger than 0.5, both SY and LC's powers started to decline quickly and were surpassed by the CLX, which were consistent with the results of Figure \ref{Fig_DGP_1_r} that the $\ell_2$-tests without regularization incorporated too many uninformative dimensions and lowered their signal to noise ratios.
We also observe that as the level of the sparsity was in the range of $[0.4,0.7]$, the extend of the power advantage of the proposed test over the other three tests became larger, which may be viewed as another confirmation of the theoretical results of the {MTT}. 

\begin{figure}
 \setlength{\abovecaptionskip}{0pt}
\setlength{\belowcaptionskip}{10pt}  \caption{Empirical powers with
respect to the signal strength $r$ for the tests of
 \cite{Srivastava_Yanagihara_2010} (SY),
\cite{Li_Chen_2012} (LC), \cite{Cai_Liu_Xia_2013} (CLX) and  the proposed multi-level thresholding test with the bootstrap calibration (MTT-BT)
for Designs 1 and 2 with Gaussian innovations under $\beta=0.6$
when $p=277$, $n_1=n_2=80$ and $p=396$, $n_1=n_2=100$ respectively.}\centering
\includegraphics[scale=0.23]{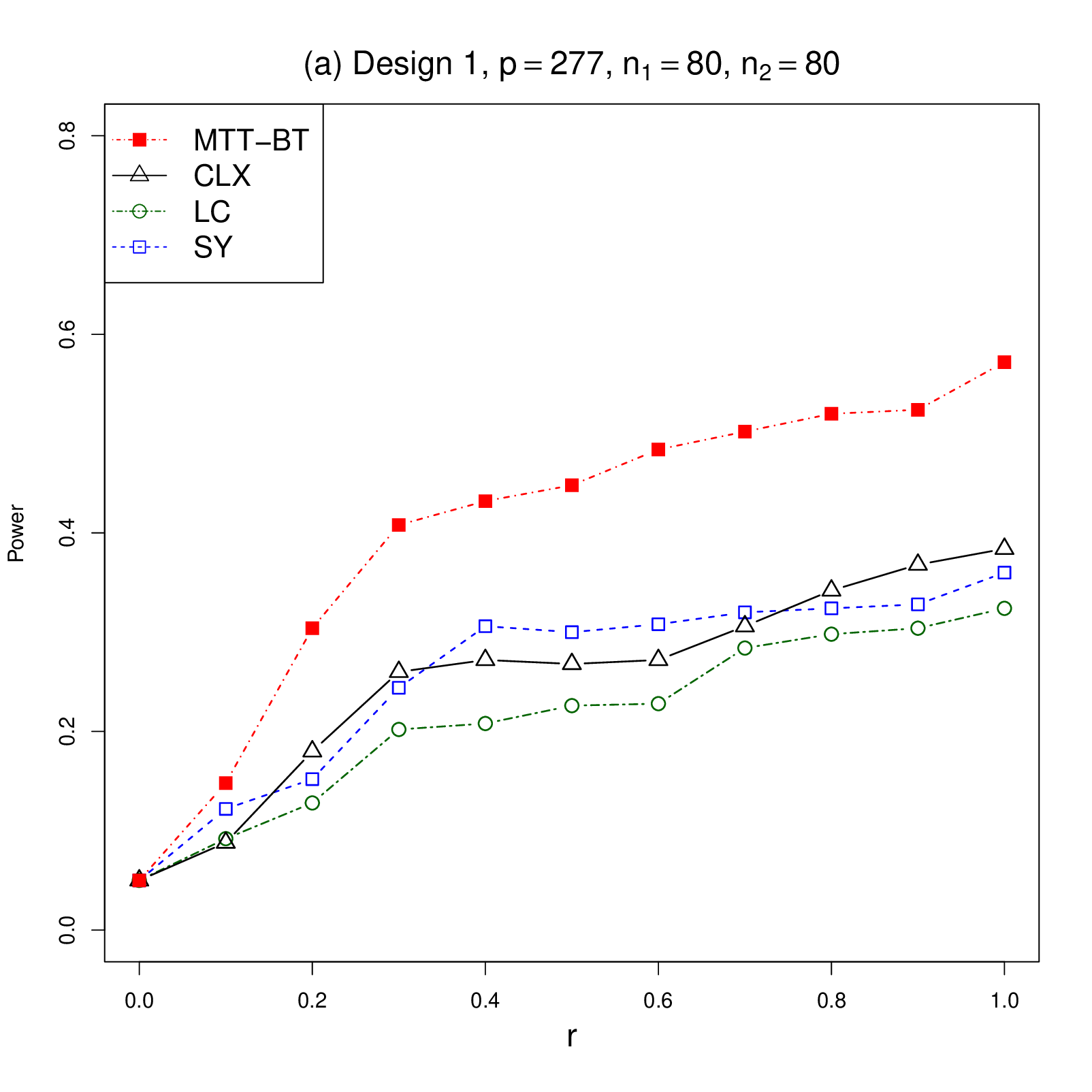}
\includegraphics[scale=0.23]{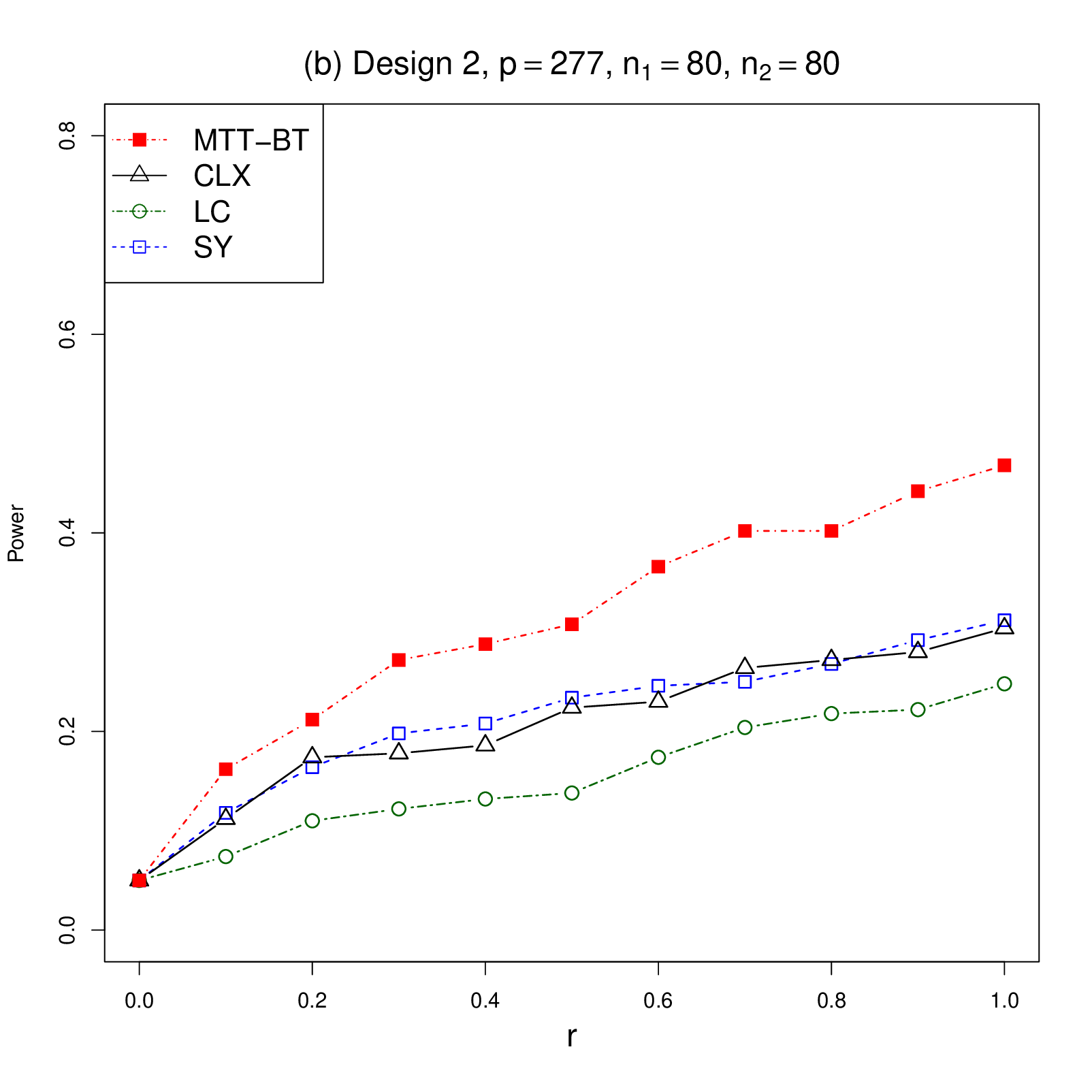}
\includegraphics[scale=0.23]{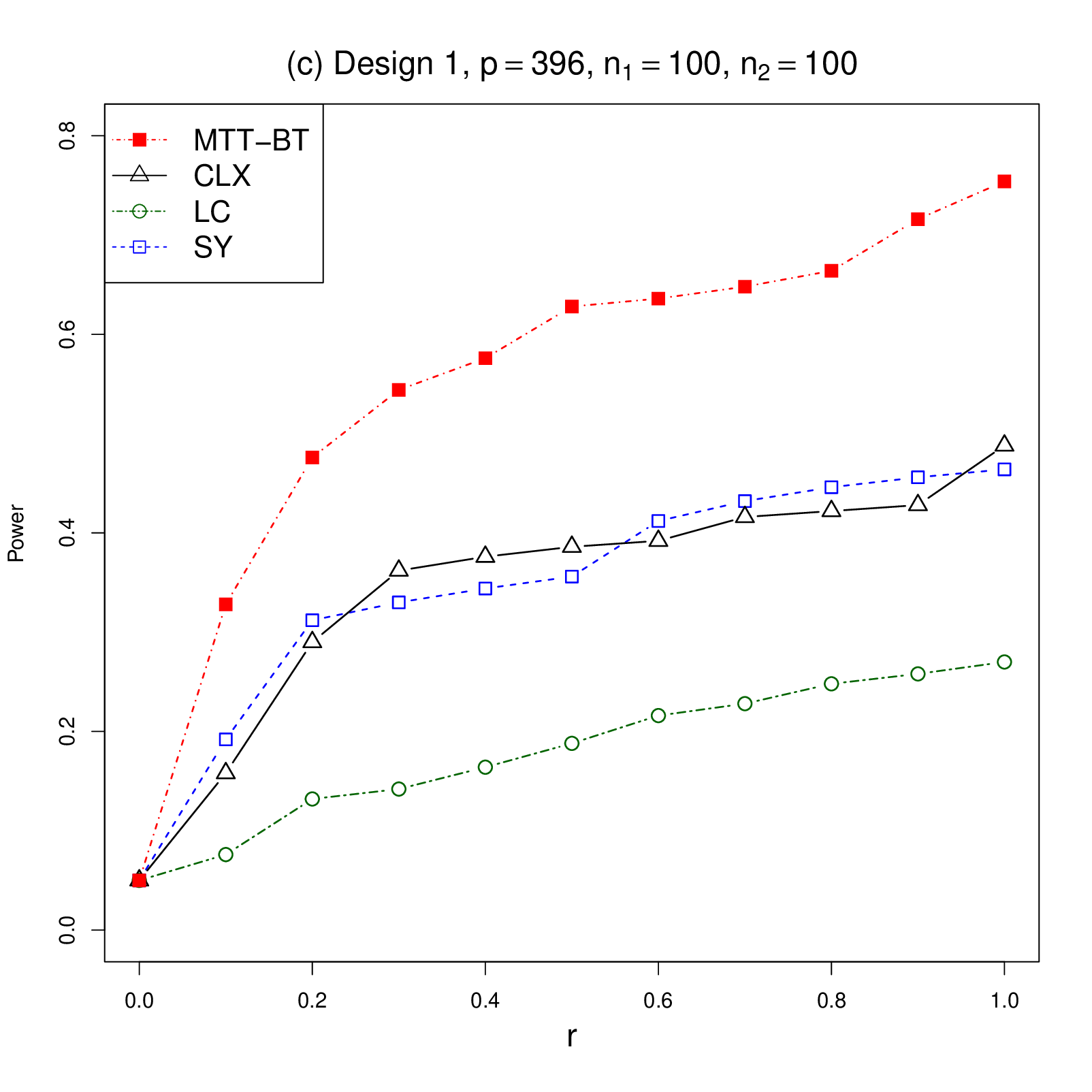}
\includegraphics[scale=0.23]{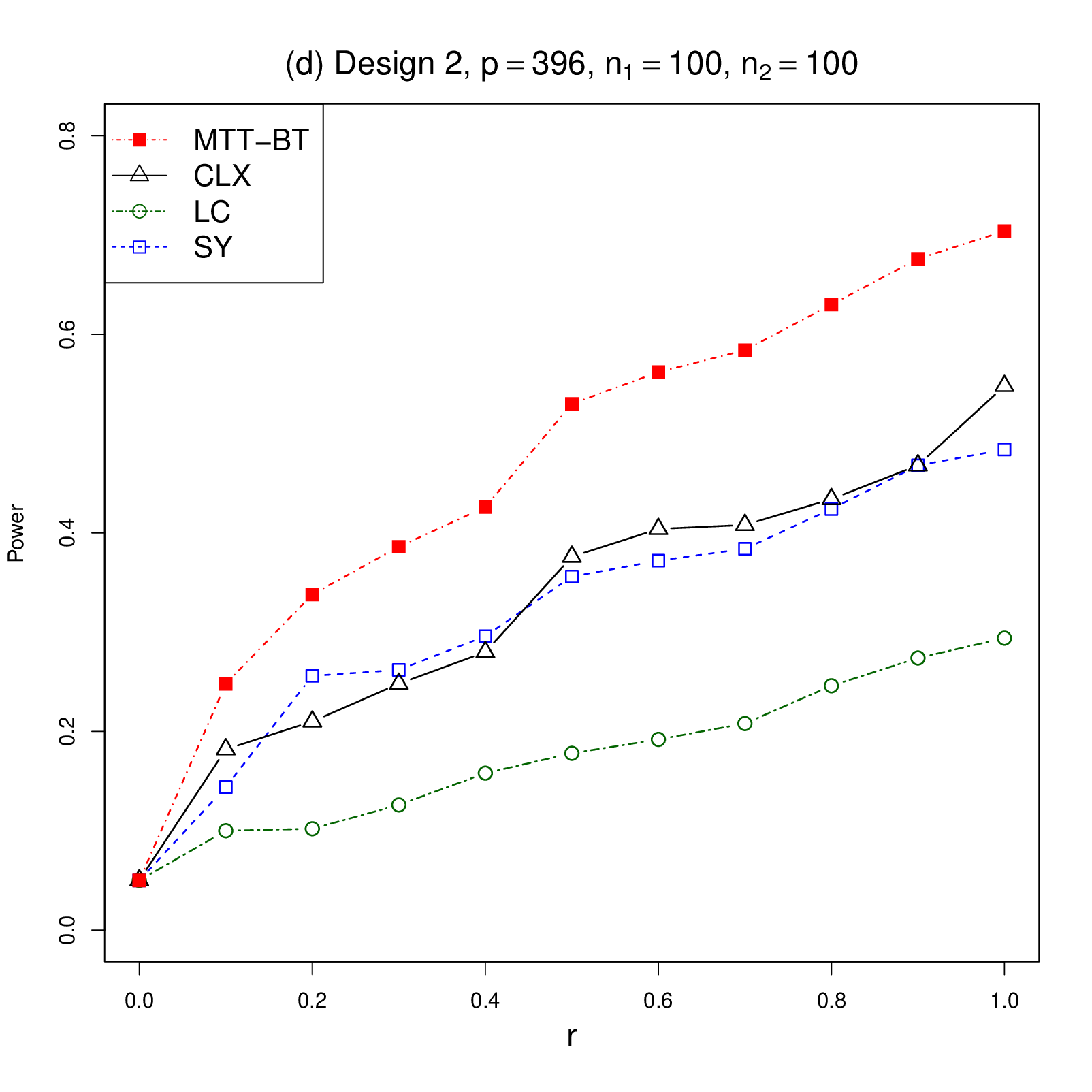}
\label{Fig_DGP_1_r}
\end{figure}

\begin{figure}
 \setlength{\abovecaptionskip}{0pt}
\setlength{\belowcaptionskip}{10pt}  \caption{Empirical powers with
respect to the sparsity level $\beta$ for the tests of
 \cite{Srivastava_Yanagihara_2010} (SY),
\cite{Li_Chen_2012} (LC), \cite{Cai_Liu_Xia_2013} (CLX) and  the proposed multi-level thresholding test with the bootstrap calibration (MTT-BT)
 for Designs 1 and 2 with Gaussian innovations under $r=0.6$
when $p=277$, $n_1=n_2=80$ and $p=396$, $n_1=n_2=100$ respectively. }\centering
\includegraphics[scale=0.23]{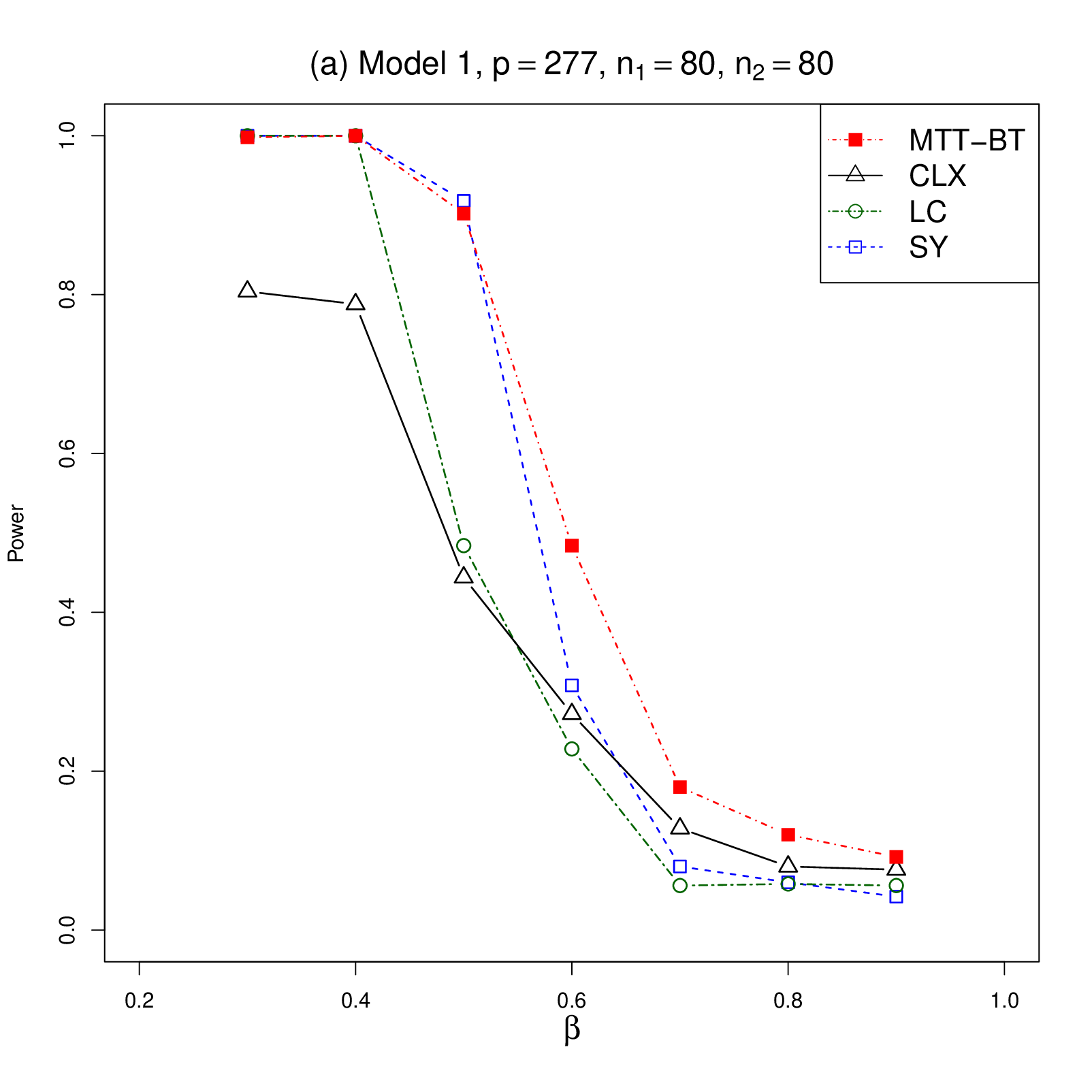}
\includegraphics[scale=0.23]{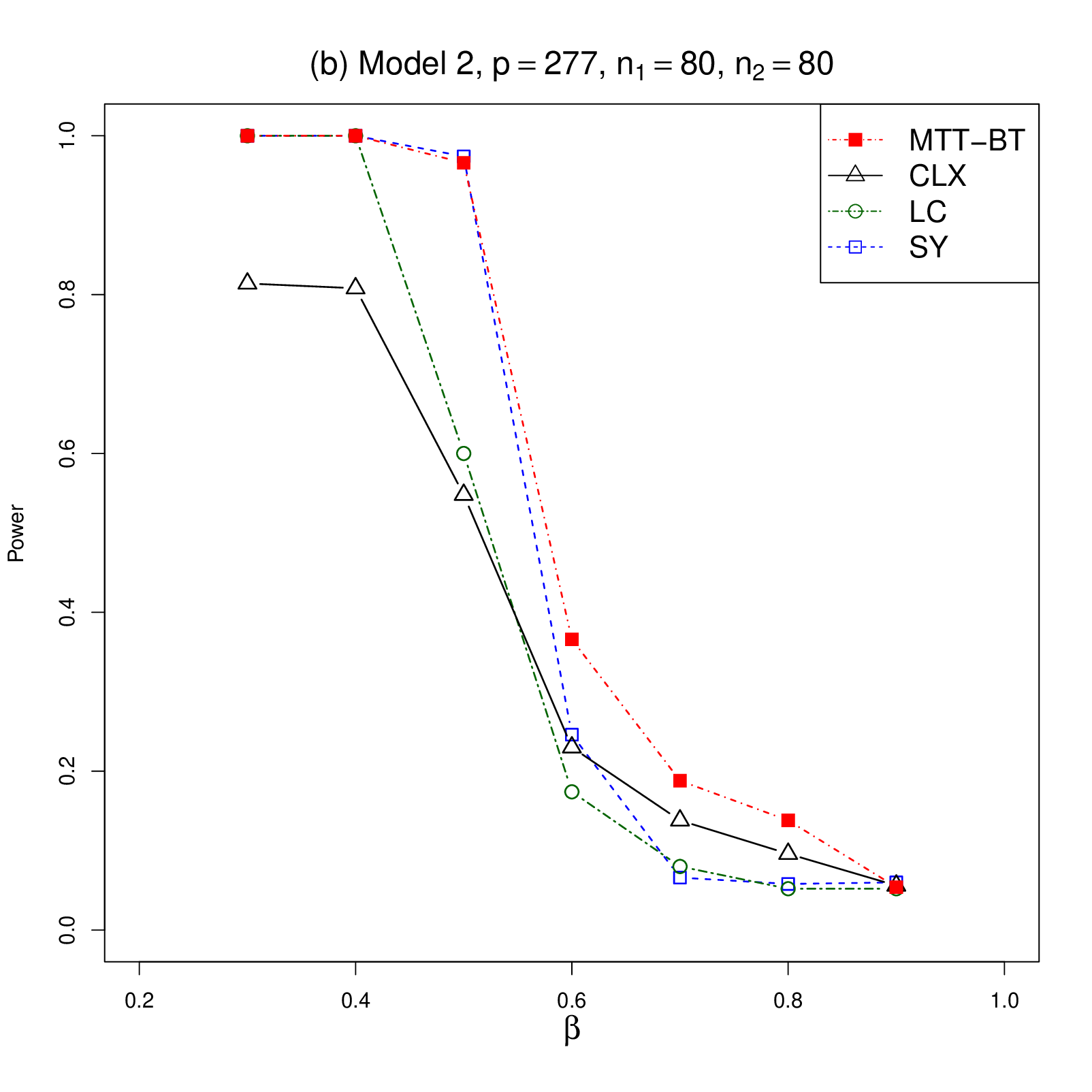}
\includegraphics[scale=0.23]{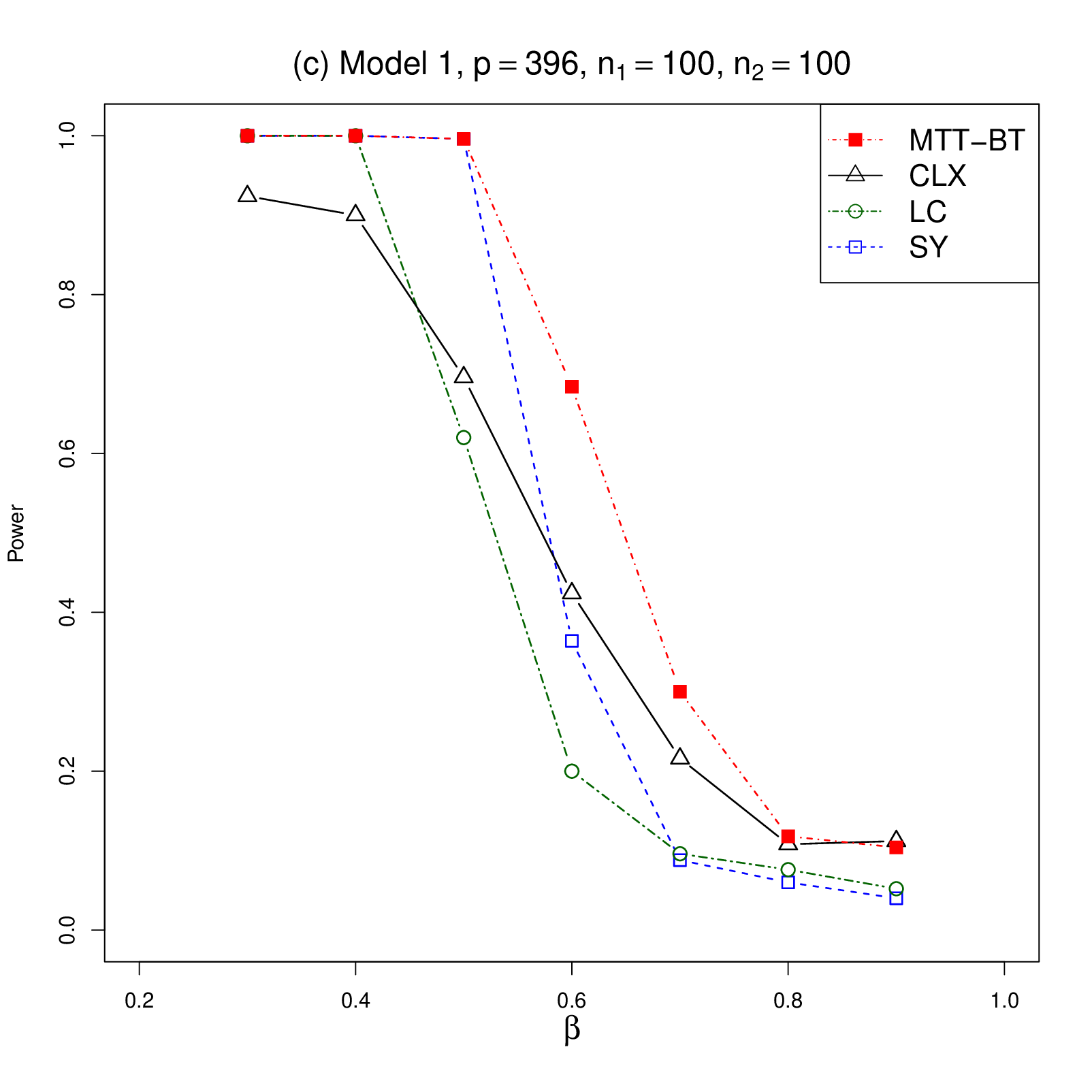}
\includegraphics[scale=0.23]{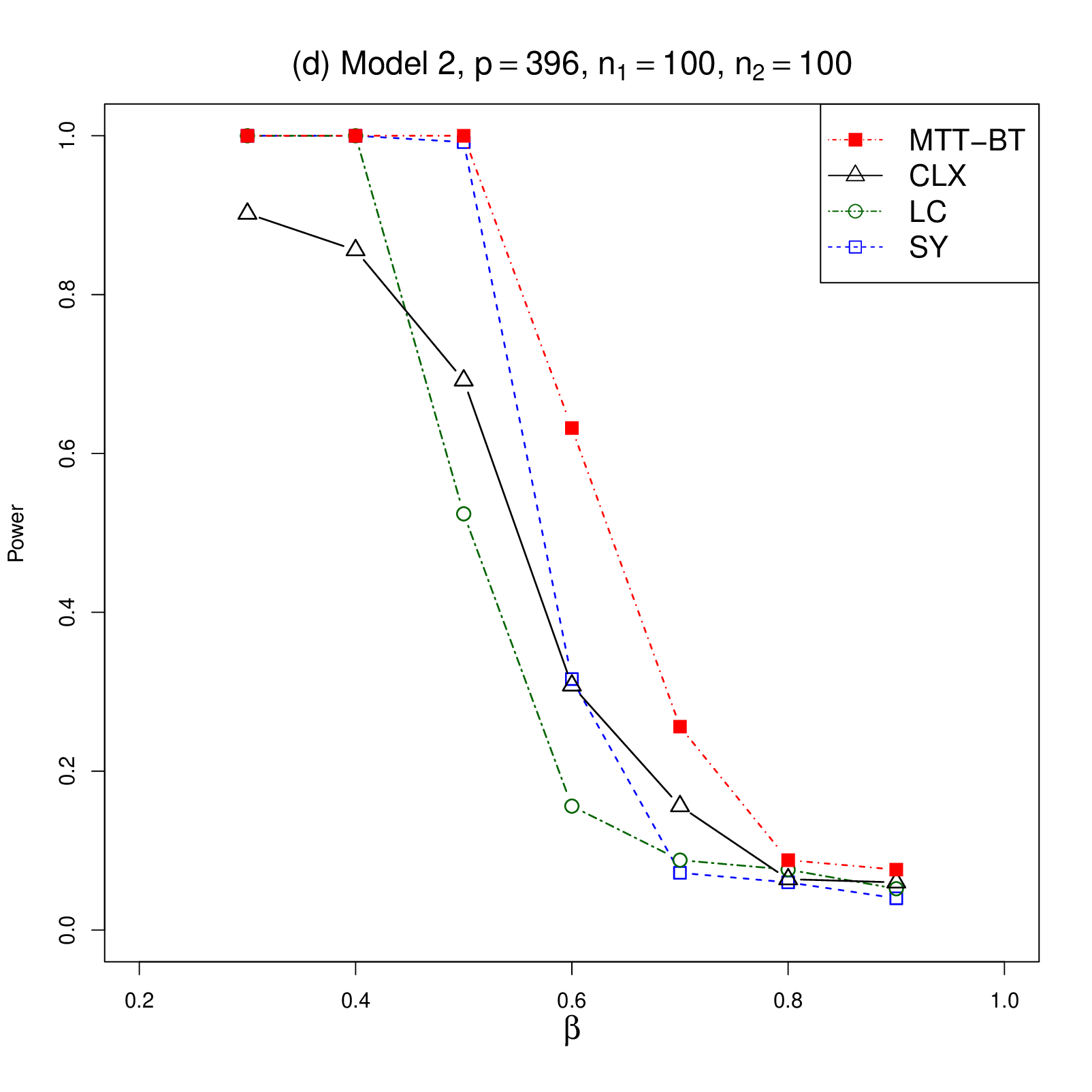}
\label{Fig_DGP_2_r}
\end{figure}

%
%
%
\setcounter{equation}{0}
\setcounter{equation}{0}

\setcounter{equation}{0}

\section{Discussion}

For establishing the asymptotic normality of the thresholding statistic $T_{n}(s)$ in (\ref{SingleTest}), the $\beta$-mixing condition (Assumption \ref{as5}) can be weakened.
Polynomial rate of the $\beta$-mixing coefficients can be assumed at the expense of more dedicated proofs.
Under this case, to prove Theorem \ref{tm1}, the length of the small and big segments of the matrix blocking ($b$ and $a$ in Figure \ref{Fig_Demo}) need to be chosen at polynomial rates of $p$, where the orders depend on the decaying rate of the $\beta$-mixing coefficients.

Although Theorem \ref{tm3} provides the minimax rate $\sqrt{\log(p) / n}$ of the signal strength for testing hypotheses (\ref{SparseH1}),  
the lower and upper bounds at this rate may  not match  due to the composite nature of the hypotheses for the two-sample test. 
To illustrate this point,
let $W \in \mathcal{W}_{\alpha}$ be the critical function of a test for the hypotheses (\ref{SparseH1}).
Let $\E_{0, \bSigma_{1}}$ and $\E_{\bSigma_{1}, \bSigma_{2}}$ be the expectation with respect to the data distribution under the null and alternative hypotheses, respectively.
The derivation of the minimax bound starts from the following 
\bea
1 + \alpha - \sup_{W \in \mathcal{W}_{\alpha}} \inf_{\bSigma_{1}, \bSigma_{2}} \E_{\bSigma_{1}, \bSigma_{2}}(W)
&\geq&
\inf_{W \in \mathcal{W}_{\alpha}} \sup_{\bSigma_{1}, \bSigma_{2}} \{ \E_{0, \bSigma_{1}} W + \E_{\bSigma_{1}, \bSigma_{2}} (1 - W) \} \nn \\
&\geq& \sup_{\bSigma_{1}} \inf_{W} \sup_{\bSigma_{2}} \{ \E_{0, \bSigma_{1}} W + \E_{\bSigma_{1}, \bSigma_{2}} (1 - W) \}. \nn
\eea
In the last inequality above, the infimum over the test $W$ is taken under fixed ${\bf\Sigma}_{1}$.
This essentially reduces to one-sample hypothesis testing.
Then, a least favorable prior will be constructed on ${\bf\Sigma}_{2}$ given ${\bf\Sigma}_{1}$ is known.
As the test under known ${\bf\Sigma}_{1}$ cannot control the type I error for the two-sample hypotheses (\ref{SparseH1}), the bound on the maxmin power for hypotheses (\ref{SparseH1}) is not tight.
It is for this reason that one may derive the tight minimax bound for the one sample spherical hypothesis $H_{0}: {\bf \Sigma} = \sigma \mathbf{I}$.
We will leave this problem as a future work, especially for the unexplored region $1/2 < \beta < \max\{2/3, (3 - \xi) / 4\}$ in Theorem \ref{tm3}.

The proposed thresholding tests can be extended to testing for correlation matrices between the two populations. Recall that $\bPsi_{1} = (\rho_{ij1})_{p\times p}$ and $\bPsi_{2} = (\rho_{ij2})_{p\times p}$ are correlation matrices of $\bX_{k}$ and $\bY_{k}$. Consider the hypotheses
$$H_0:{\bf \Psi}_1={\bf \Psi}_2  \quad \text{vs.} \quad
H_a:{\bf \Psi}_1\neq{\bf \Psi}_2.$$
Let $\hat{\rho}_{ij1} = \hat{\sigma}_{ij1} / (\hat{\sigma}_{ii1}\hat{\sigma}_{jj1})^{1/2}$ and $\hat{\rho}_{ij2} = \hat{\sigma}_{ij2} / (\hat{\sigma}_{ii2}\hat{\sigma}_{jj2})^{1/2}$ be the sample correlations of the two groups.
As for $M_{ij}$, the squared standardized difference $M_{ij}^{\star}$ between the sample correlations  can be constructed based on $\hat{\rho}_{ij1}$ and $\hat{\rho}_{ij2}$ and their estimated variances.
Let
$$T_n^{\star}(s) = \sum_{1\leq i\leq j\leq p} M_{ij}^{\star}\mathbb{I}\{M_{ij}^{\star} > \lambda_{p}(s)\}$$
be the single level thresholding statistic based on the sample correlations. Similar to the case of sample covariances, the moderate deviation results on $\hat{\rho}_{ij1} - \hat{\rho}_{ij2}$ can be derived. It can be shown that $T_n^{\star}(s)$ has the same asymptotic distribution as $T_n(s)$.
The multi-thresholding test can be constructed similar to (\ref{sup_test}) and (\ref{testproc}).

\setcounter{equation}{0}

\def\theequation{A.\arabic{equation}}
\section{Appendix}

In this section, we provide proof to Theorem \ref{tm1}. The theoretical proofs for all other propositions and theorems are relegated to the supplementary material.
Without loss of generality, we assume
$E(\bX_1)=0$ and $E(\bY_1)=0$.
Let $C$ and $L_p$ be a constant and a multi-$\log(p)$ term which may change from case to case, respectively.

\medskip

\noindent {\bf Proof of Theorem \ref{tm1}}.
To prove Theorem \ref{tm1}, we propose a novel technique that  constructs  an equivalent  U-statistic  to $T_{n}(s)$ which
is based on a partition of covariance into a group of big square blocks separated by small strips as shown in Figure \ref{Fig_Demo}.
Specifically, the indices $\{1, \ldots, p\}$ are grouped into a sequence of  big segments of length $a$ and small segments of length  $b$:
$$\{ 1, \dots, a\}, \{a+1, \dots, a + b\}, \{a + b + 1, \dots, 2a + b\}, \{2a + b + 1, \dots, 2a + 2b\}, \dots$$
where
$b = o(a)$.
Let $d = \lfloor p / (a + b) \rfloor$ be the total number of pairs of large and small segments. The sets of indices for the large segments and the small segments are, respectively,
$$S_{m} = \{(m - 1)(a + b) + 1, \dots, m a + (m - 1)b\} \mbox{ \ and \ }$$
$$ R_{m} = \{m a + (m - 1)b + 1, \dots, m(a + b)\}$$
for $m = 1, \dots, d$, and a remainder set $R_{d + 1} = \{d(a + b) + 1, \dots, p\}$.
For a two dimensional array $\{(i, j): 1 \leq i \leq j \leq p\}$, the above index partition  results in $d (d-1)/2$  square index blocks of size $a \times a$: $\{\mathcal{I}_{m_1 m_2} = S_{m_1} \times S_{m_2}: 1 \leq m_1 < m_2 \leq d\}$, colored in blue in Figure \ref{Fig_Demo}. They are separated by
$d$  smaller horizontal and vertical rectangles with  widths $a$ by $b$ and square blocks of size $b$.
There are also $d$ residual triangular blocks with $a (a + 1) / 2$ elements along the main diagonal.
The blocking scheme is 
demonstrated in Figure \ref{Fig_Demo}.

\begin{figure}[h]
\centering
\includegraphics[scale=0.5]{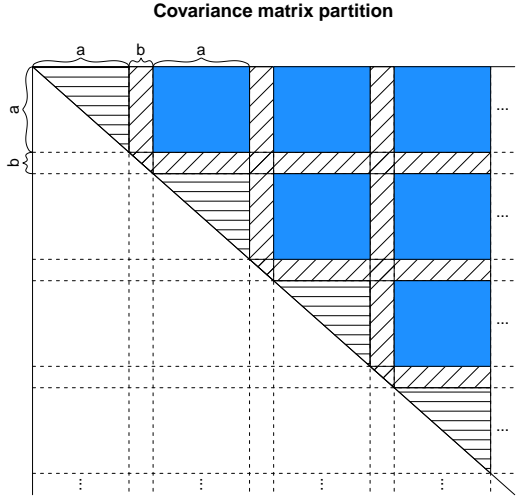}
\setlength{\abovecaptionskip}{-30pt}
\setlength{\belowcaptionskip}{0pt}
\caption{Matrix partition in the upper triangle of a covariance matrix. The square sub-matrices  (in blue color) of size $a$ are the bigger blocks, which are separated by smaller size strips of width $b$ (marked by the 45-degree lines).
There are $d$ triangle blocks along the diagonal plus remaining smaller size blocks {in the residual set $R_{d+1}$} which are not shown in the diagram. }
\label{Fig_Demo}
\end{figure}

Let
$A_{ij}(s) = L_{ij}(s) - \mu_{0,ij}(s)$, where $\mu_{0,ij}(s) = E(L_{ij}(s) | H_{0})$ and $L_{ij}(s) = M_{ij}\mathbb{I}(M_{ij} > \lambda_{p}(s))$.
Then, $T_{n}(s) - E\{T_{n}(s)\} = \sum_{1 \leq i \leq j \leq p} A_{ij}(s)$ under the null hypothesis. Here, we drop the threshold level $s$ in the notations $A_{ij}(s)$, $L_{ij}(s)$ and $\mu_{0,ij}(s)$ for simplicity of the statement, when there is no confusion.
Based on the matrix partition in Figure \ref{Fig_Demo}, $T_{n}(s) - E\{T_{n}(s)\}$ can be divided into summation of $A_{ij}(s)$ over the big square blocks of size $a \times a$, the small strips and the triangular blocks along the main diagonal.

Let $R = \cup_{m=1}^{d} R_{m}$ be the collection of the indices in the small segments.
From Figure \ref{Fig_Demo}, $T_{n}(s) - E\{T_{n}(s)\}$ can be divided into four parts such that
\be T_{n}(s) - E\{T_{n}(s)\} = B_{1,n} + B_{2,n} + B_{3,n} + B_{4,n}, \label{eq:decom}
\ee
where
\bea
B_{1, n} = \sum_{1\leq m_1<m_2 \leq d} \ \sum_{i \in S_{m_1},~ j \in S_{m_2}} A_{ij}, &&
B_{2, n} = \sum_{i \in R ~ or ~ j \in R, ~i \leq j} A_{ij}, \nn \\
B_{3, n} = \sum_{1 \leq m \leq d} \ \sum_{i,j \in S_{m}, ~i \leq j} A_{ij}, &&
B_{4, n} = \sum_{j \in R_{d+1}, ~i \leq j} A_{ij}, 
\label{eq:Decomposition}\eea
Here,  $B_{1, n}$ is the sum of  $A_{ij}$ over the $d(d-1)/2$ big $a \times a$ square blocks in Figure \ref{Fig_Demo},  $B_{2, n}$ is the sum over all the smaller rectangular and square blocks, $B_{3, n}$ is over the $d$ triangular blocks along the main diagonal,
and $B_{4, n}$ is over the remaining segments  including the residual blocks $R_{d + 1}$ towards the right end columns  of the matrix,  
 respectively.

For the decomposition of $T_{n}(s) - E\{T_{n}(s)\}$ in (\ref{eq:Decomposition}),
let $N_{l}$ be the number of elements in $B_{l, n}$ for $l = 1, \dots, 4$.
Note that $N_{1} = a^2 d (d - 1) / 2 = q(1 + o(1))$,
$N_2 \leq d p b \leq p^2 b / (a + b) = o(p^{2})$, $N_{3} = d a^2 / 2 \leq p a / 2 = o(p^2)$ and $N_{4} \leq |R_{d+1}| p \leq (a + b)p = o(p^2)$. 
Similar as deriving $\text{Var}\{T_{n}(s)\}$ in Proposition \ref{pn1}, we have
\bea
\text{Var}(B_{2, n}) &=&
\sum_{i \in R ~or~ j \in R, ~i \leq j} \text{Var}(A_{ij}) \label{eq:B2Var} \\
&+& \mbox{Cov}\bigg( \sum_{\substack{i_{1} \in R ~or~ j_{1} \in R \\ i_{1} \leq j_{1}}} A_{i_1j_1},
\sum_{\substack{i_{2} \in R ~or~ j_{2} \in R \\ i_{2} \leq j_{2}}} A_{i_2j_2} \bigg), \label{eq:B2Cov}
\eea
where $\mbox{Var}(A_{ij}) = \mbox{Var}(L_{ij}) = v(0, s)\{1 + o(1)\}  \sim L_{p}p^{-2s}$ under the null hypothesis, which is given in Lemma 5 in the SM. Notice that the summation of the variances on the right side of (\ref{eq:B2Var}) is bounded by $L_{p}p^{-2s} N_{2} = o\big(\text{Var}\{T_{n}(s)\}\big)$.

For the covariance terms in (\ref{eq:B2Cov}),
let $d_{i_1j_1,i_2j_2} = \min(|i_1 -i_2|, |i_1 -j_2|, |j_1-j_2|, |j_1-i_2|)$ be the minimum coordinate distance between $(i_1,j_1)$ and $(i_2,j_2)$, and between $(i_1,j_1)$ and $(j_2,i_2)$, where $i_{1} \leq j_{1}$ and $i_{2} \leq j_{2}$.
For any fixed $(i_1,j_1)$ and a large positive constant $M$, by Assumption \ref{as5} and Davydov's inequality (Corollary 1.1 of \cite{Bosq_1998}, p.21), there exists a constant $c > 0$ such that $|\mbox{Cov}(L_{i_1j_1}, L_{i_2j_2})| \leq C\gamma_{1}^{d_{i_1j_1,i_2j_2}} \leq p^{-M}$ for $\gamma_{1} \in (0, 1)$ and any $d_{i_1j_1,i_2j_2} > c \log p$.
Therefore,
\bea
&&\bigg| \mbox{Cov}\bigg( \sum_{i_{1} \in R ~or~ j_{1} \in R, ~i_{1} \leq j_{1}} A_{i_1j_1}, \sum_{i_{2} \in R ~or~ j_{2} \in R, ~i_{2} \leq j_{2}} A_{i_2j_2} \bigg) \bigg| \nn \\
&\leq&  \sum_{i_{1} \in R ~or~ j_{1} \in R, ~i_{1} \leq j_{1}} \bigg\{ N_{2}p^{-M}
 + \sum_{\substack{d_{i_1j_1,i_2j_2} \leq c \log p \\ i_{2} \in R ~or~ j_{2} \in R, ~i_{2} \leq j_{2}}} \big| \mbox{Cov}(A_{i_1j_1},  A_{i_2j_2} )\big| \bigg\}, \nn
\eea
where by Lemmas 5 and 6 in the SM, $|\mbox{Cov}(A_{i_1j_1},  A_{i_2j_2} )\big| \leq L_{p}|\rho_{i_1j_1,i_2j_2}| p^{-\frac{4s}{2 - \epsilon}} + \mu_{0,i_1j_1}\mu_{0,i_2j_2}L_{p}n^{-1/2}$ for a small $\epsilon > 0$.
It follows that
\bea
\sum_{\substack{d_{i_1j_1,i_2j_2} \leq c \log p \\ i_{2} \in R ~or~ j_{2} \in R, ~i_{2} \leq j_{2}}} \big| \mbox{Cov}(A_{i_1j_1},  A_{i_2j_2} )\big|
&\leq& \sum_{|i_2 - i_1| \leq c \log p}\sum_{j_2 = 1}^{p} \big| \mbox{Cov}(A_{i_1j_1},  A_{i_2j_2} )\big| \nn \\
&+& \sum_{|j_2 - j_1| \leq c\log p}\sum_{i_2 = 1}^{p} \big| \mbox{Cov}(A_{i_1j_1},  A_{i_2j_2} )\big|, \nn
\eea
which is bounded by $L_{p}\sum_{j_2 = 1}^{p} |\rho_{i_1j_1,i_2j_2}|p^{-\frac{4s}{2 - \epsilon}} + L_{p}p^{1-4s}n^{-1/2}$.
It has been shown that $\sum_{j_2 = 1}^{p} |\rho_{i_1j_1,i_2j_2}| \leq C < \infty$ in (S.23) in the SM.
By choosing $M$ large, the covariance term in (\ref{eq:B2Cov}) is bounded by $L_{p}N_{2}p^{-\frac{4s}{2 - \epsilon}} + L_{p}N_{2}p^{1-4s}n^{-1/2}$, which is a small order term of $\text{Var}\{T_{n}(s)\}$ if $s > 1/2$ under Assumption \ref{as1} or $s > 1/2 - \xi / 4$ under Assumption \ref{as1poly}. Therefore, $\text{Var}(B_{2, n})= o\big(\text{Var}\{T_{n}(s)\}\big)$.


For $B_{3, n}$, note that the triangles $\{(i, j)\in S_{m}, i \leq j\}$ along the diagonal are at least $b$ apart from each other, where $b \sim \log(p)$. The covariances between $\sum_{i,j \in S_{m_{1}}} A_{ij}$ and $\sum_{i,j \in S_{m_{2}}} A_{ij}$ are negligible for $m_1 \neq m_2$. It follows that
\bea
\mbox{Var}(B_{3, n}) &=& \sum_{1 \leq m \leq d} \mbox{Var}\bigg(\sum_{i,j \in S_{m}, ~i \leq j} A_{ij}\bigg) \{1 + o(1)\} \nn \\
&=& \sum_{1 \leq m \leq d}\sum_{i_{1},j_{1} \in S_{m}, ~i_{1} \leq j_{1}}\sum_{i_{2},j_{2} \in S_{m}, ~i_{2} \leq j_{2}} \mbox{Cov}(A_{i_{1}j_{1}}, A_{i_{2}j_{2}}) \{1 + o(1)\}, \nn
\eea
which is bounded by $C d a^{4} v(0, s) = o( L_{p} a^{3} p^{1-2s})$. This shows that $\mbox{Var}(B_{3, n}) = o\big(\text{Var}\{T_{n}(s)\}\big)$ when $a \ll p^{1/3}$.
Here, for two positive sequences $\{c_{1, n}\}$ and $\{c_{2, n}\}$, $c_{1, n} \ll c_{2, n}$ means that $c_{1, n} = o(c_{2, n})$.
For $B_{4, n}$, we have
\bea
\mbox{Var}(B_{4, n}) 
&\leq& N_{4}v(0, s)\{1 + o(1)\} + \sum_{j_{1} \in R_{d+1}}\sum_{j_{2} \in R_{d+1}}|\mbox{Cov}(A_{i_{1}j_{1}}, A_{i_{2}j_{2}})| \nn \\
&=& o(p^{2 - 2s}) + \sum_{j_{1} \in R_{d+1}}\bigg( N_{4}p^{-M} + \sum_{\substack{d_{i_1j_1,i_2j_2} \leq c \log p \\ j_{2} \in R_{d+1}}}|\mbox{Cov}(A_{i_{1}j_{1}}, A_{i_{2}j_{2}})|\bigg). \nn
\eea
Similar to the case of $\mbox{Var}(B_{2, n})$, the last summation term above is equal to $N_{4}^{2}p^{-M} + L_{p}N_{4}p^{-\frac{4s}{2 - \epsilon}} + L_{p}N_{4}p^{1-4s}n^{-1/2}$, which is a small order term of $\text{Var}\{T_{n}(s)\}$.
Meanwhile, since $N_{1} = q(1 + o(1))$, following the same derivation of Proposition \ref{pn1}, it can be shown that $\text{Var}(B_{1, n}) = \text{Var}(T_{n}(s))\{1 + o(1)\}$.

Combining the above results, we see that $\text{Var}(B_{l, n})$ are at a smaller order of $\text{Var}\{T_{n}(s)\}$ for $l = 2, \ldots, 4$.
This together with (\ref{eq:decom}) imply
\be\frac{T_{n}(s) - E(T_{n}(s))}{\sqrt{\text{Var}(T_{n}(s))}} = \frac{B_{1, n}}{\sqrt{\text{Var}(T_{n}(s))}} + o(1).\label{eq:BlockMainOrder}\ee
Therefore, to show the asymptotical normality of $T_{n}(s) - E\{T_{n}(s)\}$, it suffices to focus on its main order term $B_{1, n}$.

Let ${\bf Z}_{S_{m}} = \{{\bf X}_{S_{m}}, {\bf Y}_{S_{m}}\}$ for $m = 1, \ldots, d$, where ${\bf X}_{S_{m}} = \{X_{ki}: 1 \leq k \leq n_1, i \in S_{m}\}$ and ${\bf Y}_{S_{m}} = \{Y_{ki}: 1 \leq k \leq n_2, i \in S_{m}\}$ are the segments of the two data matrices with the columns in $S_{m}$.
Notice that the summation of $A_{ij}$ in $B_{1, n}$ 
can be expressed as
$$\sum_{i \in S_{m_1}, j \in S_{m_2}} A_{ij} = f({\bf Z}_{S_{m_1}}, {\bf Z}_{S_{m_2}})$$
for some function $f(\cdot, \cdot)$.

Let $\mathcal{F}_{m_a}^{m_b}({\bf Z})=\sigma\{{\bf Z}_{S_{m}}: m_a \leq m \leq m_b\}$ be the $\sigma$-algebra generated by $\{{\bf Z}_{S_{m}}\}$ for $1 \leq m_a \leq m_b \leq d$.
Let $\zeta_{z}(h) = \sup_{1 \leq m \leq d - h} \zeta\{\mathcal{F}_{1}^{m}({\bf Z}), \mathcal{F}_{m+h}^{d}({\bf Z})\}$ be the $\beta$-mixing coefficient of the sequence ${\bf Z}_{S_{1}}, \dots, {\bf Z}_{S_{d}}$.
By Theorem 5.1 in \cite{Bradley_2005} and Assumption \ref{as5}, we have
$$\zeta_{z}(h) \leq \sum_{k_1=1}^{n_1}\zeta_{x, p}(hb) + \sum_{k_2=1}^{n_2}\zeta_{y, p}(hb) \leq C(n_1 + n_2)\gamma^{hb}$$
for some $\gamma\in(0,1)$. Choosing $b = b_{0} - \log(n_1 + n_2) / \log(\gamma)$ leads to $\zeta_{z}(h) \leq C \gamma^{hb_{0}}(n_1 + n_2)^{1-h} \leq C \gamma^{hb_{0}}$.
By Berbee's theorem (page 516 in \cite{AthreyaLahiri}), there exist ${\bf Z}^{^{\ast}}_{S_{2}}$ independent of ${\bf Z}_{S_{1}}$ such that $P({\bf Z}_{S_{2}} \neq {\bf Z}^{^{\ast}}_{S_{2}}) = \zeta\{\sigma({\bf Z}_{S_{1}}), \sigma({\bf Z}_{S_{2}})\} \leq \zeta_{z}(1) \leq C\gamma^{b_{0}}$.
By applying this theorem recursively, there exist ${\bf Z}_{S_{1}}, {\bf Z}^{^{\ast}}_{S_{2}}, \dots, {\bf Z}^{^{\ast}}_{S_{d}}$ that are mutually independent with each other, and $P({\bf Z}_{S_{m}} \neq {\bf Z}^{^{\ast}}_{S_{m}}) \leq C\gamma^{b_{0}}$ for $m = 2, \dots, d$.

Let $D = \cup_{m=2}^{d}\{{\bf Z}_{S_{m}} \neq {\bf Z}^{^{\ast}}_{S_{m}}\}$, then $P(D) \leq C d \gamma^{b_{0}}$.
By choosing $b_{0} = -c_{0} \log(p) / \log(\gamma)$ for a large positive number $c_{0}$, we have $P(D)$ converges to $0$ at the rate $p^{1 - c_{0}} / (a + b)$.
Notice that $\text{Var}(B_{1, n})$ is at the order $p^{2 - 2s}$.
Since $E(|B_{1, n} \mathbb{I}_{D}|) \leq \{\text{Var}(B_{1, n}) P(D)\}^{1/2}$ converges to 0 for a sufficiently large $c_{0}$,
it follows that $B_{1, n} \mathbb{I}_{D} \to 0$ in probability by choosing a large $c_{0}$.

Thus, by letting $b = c_{1}\max\{\log(p), \log(n_1+n_2)\}$ for a large constant $c_{1}>0$, there exists an array of mutually independent random vectors ${\bf Z}_{S_{1}}^{\ast}, \dots, {\bf Z}_{S_{d}}^{\ast}$ such that ${\bf Z}_{S_{m}}^{\ast} = {\bf Z}_{S_{m}}$ with overwhelming probability for $m = 1, \ldots, d$ and $B_{1, n}$ can be expressed as a $U$-statistic formulation on a sequence of mutually independent random vectors as
\be
B_{1, n} = \sum_{m_1 < m_2} f({\bf Z}_{S_{m_1}}^{\ast},{\bf Z}_{S_{m_2}}^{\ast}).
\label{Ustat}\ee
For simplicity of notations, we will drop the superscript $\ast$ in (\ref{Ustat}) in the following proof.
Now, we only need to establish the asymptotical normality of $B_{1, n}$ under the expression (\ref{Ustat}).

To this end, we first study the conditional distribution of $M_{ij}$ given the $j$th variable, where $i \in S_{m_1}$, $j \in S_{m_2}$ and $m_1 \neq m_2$.
Recall that $$F_{ij}=(\hat{\sigma}_{ij1}-\hat{\sigma}_{ij2})(\hat{\theta}_{ij1}/n_1+\hat{\theta}_{ij2}/n_2)^{-1/2}$$ is the standardization of $\hat{\sigma}_{ij1}-\hat{\sigma}_{ij2}$, where $\hat{\sigma}_{ij1}=\tilde{\sigma}_{ij1}-\bar{X}_i\bar{X}_j$ and $\hat{\sigma}_{ij2} = \tilde{\sigma}_{ij2}-\bar{Y}_i\bar{Y}_j$
for $\tilde{\sigma}_{ij1}=\sum_{k=1}^{n_1}X_{ki}X_{kj} / n_1$ and $\tilde{\sigma}_{ij2}=\sum_{k=1}^{n_2}Y_{ki}Y_{kj} / n_2$. Then, $M_{ij} = F_{ij}^{2}$.
Note that the unconditional asymptotical distribution of $F_{ij}$ is standard normal.

Let $E_{j}(\cdot)$, $\text{Var}_{j}(\cdot)$ and $\text{Cov}_{j}(\cdot)$ be the conditional mean, variance and covariance given the $j$th variable, respectively.
From the proof of Lemma 7 in the SM,
we have that $E_{j}(\hat{\sigma}_{ij1}) = E_{j}(\hat{\sigma}_{ij2}) = 0$, $\text{Var}_{j}(\tilde{\sigma}_{ij1}) = \sigma_{ii1}\tilde{\sigma}_{jj1} / n_1$ and $\text{Cov}_{j}(\tilde{\sigma}_{ij1}, \bar{X}_i\bar{X}_j) = \text{Var}_{j}(\bar{X}_i\bar{X}_j) = \sigma_{ii1}(\bar{X}_{j})^{2} / n_1$.
It follows that $\text{Var}_{j}(\hat{\sigma}_{ij1}) = \sigma_{ii1}\hat{\sigma}_{jj1} / n_1$ and $\text{Var}_{j}(\hat{\sigma}_{ij2}) = \sigma_{ii2}\hat{\sigma}_{jj2} / n_2$.
In the proof of Lemma 7, it has also been shown that $\hat{\theta}_{ij1} = \sigma_{ii1}\hat{\sigma}_{jj1} + O_{p}(\sqrt{\log(p) / n})$ and $\hat{\theta}_{ij2} = \sigma_{ii2}\hat{\sigma}_{jj2} + O_{p}(\sqrt{\log(p) / n})$ given the $j$th variable.
Similar results hold given the $i$th variable.
Therefore, $F_{ij}$ is still asymptotically standard normal distributed given either the $i$th or the $j$th variable.
And, the moderate deviation results from Lemma 2.3 and Theorem 3.1 in \cite{Saulis_1991} for independent but non-identically distributed variables can be applied to $F_{ij}$, given either one of the variables.

Let $\mathscr{F}_0=\{\emptyset,\Omega\}$ and $\mathscr{F}_m=\sigma\{{\bf Z}_{S_{1}},\dots,{\bf Z}_{S_{m}}\}$ for $m=1,2,\cdots,d$ be a sequence of $\sigma$-field generated by $\{{\bf Z}_{S_{1}},\dots,{\bf Z}_{S_{m}}\}$. Let $E_{\mathscr{F}_m}(\cdot)$ denote the conditional expectation with respect to $\mathscr{F}_m$.
Write $B_{1, n}=\sum^d_{m=1}D_{m}$, where $D_{m}=(E_{\mathscr{F}_m} - E_{\mathscr{F}_{m-1}})B_{1, n}$.
Then for every $n, p$, $\{D_{m},1\leq m\leq d\}$ is a martingale difference sequence with respect to the $\sigma$-fields $\{\mathscr{F}_m\}_{m=0}^{\infty}$.
Let $\sigma^2_{m}=E_{\mathscr{F}_{m-1}}(D^2_{m})$.
By the martingale central limit theorem \citep[Chapter 3 in][]{Hall_Heyde_1980}, 
 to show the asymptotical normality of $B_{1, n}$, it suffices to show
\be
\frac{\sum^d_{m=1}\sigma^2_{m}}{\text{Var}(B_{1, n})}\xrightarrow[]{p}1\quad\mbox{ \ and \ }\quad\frac{\sum^d_{m=1} E(D^4_{m})}{\text{Var}^{2}(B_{1, n})}\xrightarrow[]{}0.
\label{MCLT}\ee

By the independence between $\{{\bf Z}_{S_{1}},\dots,{\bf Z}_{S_{d}}\}$, we have
\bea
D_{m} &=&
\sum_{m_1 = 1}^{m-1}f({\bf Z}_{S_{m_1}},{\bf Z}_{S_{m}}) \ + \sum_{m_2>m} E_{\mathscr{F}_m} f({\bf Z}_{S_{m}},{\bf Z}_{S_{m_2}}) \label{eq:MDiff1} \\
&-& \sum_{m_1 = 1}^{m - 1}E_{\mathscr{F}_{m-1}} f({\bf Z}_{S_{m_1}},{\bf Z}_{S_{m}}),\nn
\eea
where for any $m_1 < m_2$,
$$E_{\mathscr{F}_{m_1}} f({\bf Z}_{S_{m_1}},{\bf Z}_{S_{m_2}}) = E_{\mathscr{F}_{m_1}} \sum_{i \in S_{m_1}, j \in S_{m_2}} A_{ij} = \sum_{i \in S_{m_1}, j \in S_{m_2}} E_{i} A_{ij}.$$
For $m_1 < m_2$, let
\be\begin{split}
&\tilde{f}({\bf Z}_{S_{m_1}}, {\bf Z}_{S_{m_2}}) = \sum_{i \in S_{m_1}, j \in S_{m_2}} \tilde{A}_{ij} \mbox{ \ for} \\
\tilde{A}_{ij} = M_{ij}&\mathbb{I}(M_{ij} > \lambda_{p}(s)) - E_{\mathscr{F}_{m_1}}\{M_{ij}\mathbb{I}(M_{ij} > \lambda_{p}(s))\},
\end{split}\nn\ee
where $i \in S_{m_1}$ and $j \in S_{m_2}$. We can decompose $D_{m} = D_{m, 1} + D_{m, 2}$, where
\be
D_{m, 1} = \sum_{m_1 = 1}^{m-1}\tilde{f}({\bf Z}_{S_{m_1}},{\bf Z}_{S_{m}}) \mbox{ \ and \ } D_{m, 2} = \sum_{m_2>m} E_{\mathscr{F}_m} f({\bf Z}_{S_{m}},{\bf Z}_{S_{m_2}}).
\label{eq:MDiff2}\ee

Let $G_{0} = \{\max_{k_{1}, k_{2}, i}(|X_{k_{1}i}|, |Y_{k_{2}i}|) \leq c \sqrt{\log p}\}$ for a positive constant $c$. Under Assumption \ref{as3}, $P(G_{0}^{c})\to 0$ in a polynomial rate of $p$ for a large $c$. 
To study $\{\sigma_{m}^{2}\}$, we focus on the set $G_{0}$. 
By Lemma 7 in the SM, 
we have $$E_{i} \{M_{ij}\mathbb{I}(M_{ij} > \lambda_{p}(s)) \} = \mu_{0, ij}\big\{1 + O(L_pn^{-1/2})\big\},$$
which implies $E_{i} A_{ij} = \mu_{0, ij} O(L_pn^{-1/2})$. This leads to $E_{\mathscr{F}_m} f({\bf Z}_{S_{m}},{\bf Z}_{S_{m_2}}) = L_p O(a^{2} p^{-2s} n^{-1/2})$ for $m_2 > m$,
and $D_{m, 2} = (d - m) L_p O(a^{2} p^{-2s} n^{-1/2})$.
From (\ref{eq:MDiff2}), we can write $D_{m}^{2} = D_{m, 1}^{2} + 2 D_{m, 1} D_{m, 2} + D_{m, 2}^{2}$, where $D_{m, 2}^{2} = (d - m)^{2} L_p O(a^{4} p^{-4s} n^{-1})$.
Note that $E_{\mathscr{F}_{m-1}}(D_{m, 1} D_{m, 2})$ is equal to
\bea
\sum_{m_1 = 1}^{m-1}\sum_{m_2 > m}\sum_{j_1 \in S_{m_1},j_2\in S_{m}}\sum_{j_3 \in S_{m},j_4\in S_{m_2}} E_{\mathscr{F}_{m-1}}\{\tilde{A}_{j_1j_2} E_{\mathscr{F}_{m}}(A_{j_3j_4})\}. \nn
\eea
Similar as applying the coupling method on the big segments ${\bf Z}_{S_{1}}, \dots, {\bf Z}_{S_{d}}$ of the variables,
the $j_2$th and $j_3$th variables can be effectively viewed as independent when $|j_2 - j_3| > c \log p$ for some constant $c > 0$.
Therefore, given $\mathscr{F}_{m-1}$, $E_{\mathscr{F}_{m-1}}\{\tilde{A}_{j_1j_2} E_{\mathscr{F}_{m}}(A_{j_3j_4})\}$ is negligible when $|j_2 - j_3| > c \log p$.
Meanwhile, notice that $\big| E_{\mathscr{F}_{m-1}}\{\tilde{A}_{j_1j_2} E_{\mathscr{F}_{m}}(A_{j_3j_4})\} \big| \leq  O(L_p p^{-2s} n^{-1/2}) E_{\mathscr{F}_{m-1}}(|\tilde{A}_{j_1j_2}|)$ and $E_{\mathscr{F}_{m-1}}(|\tilde{A}_{j_1j_2}|) \leq 2 E_{\mathscr{F}_{m-1}}\{M_{ij}\mathbb{I}(M_{ij} > \lambda_{p}(s))\}$, which is at the order $L_{p}p^{-2s}$.
Therefore, we have $$\big| E_{\mathscr{F}_{m-1}}(D_{m, 1} D_{m, 2}) \big| \leq O(L_{p} d^{2} a^{3} p^{-4s} n^{-1/2}).$$
%
Base on the above results, by choosing $a \ll \sqrt{n}$, $\sigma^2_{m} = E_{\mathscr{F}_{m-1}}(D^2_{m})$ can be expressed as $\sigma^2_{m} = E_{\mathscr{F}_{m-1}}(D_{m, 1}^{2}) + O(L_p d^{2} a^{3} p^{-4s} n^{-1/2})$, where
\be
E_{\mathscr{F}_{m-1}}(D_{m, 1}^{2}) = \sum_{m_1,m_2 = 1}^{m-1} \sum_{\substack{j_1\in S_{m_1} \\ j_2 \in S_{m}}} \sum_{\substack{j_3\in S_{m_2} \\ j_4 \in S_{m}}} E_{\mathscr{F}_{m-1}}(\tilde{A}_{j_1j_2}\tilde{A}_{j_3j_4}).\label{eq:MDiffVar}\ee

For the above summation in (\ref{eq:MDiffVar}), note that when $j_1 = j_3$, $j_2 = j_4$, $$E_{\mathscr{F}_{m-1}}(\tilde{A}_{j_1j_2}^2) = 
E_{j_1}\{M_{j_1j_2}^{2}\mathbb{I}(M_{j_1j_2} > \lambda_{p}(s))\} - \mu_{0, j_1j_2}^{2}(1 + o_p(1)).$$
By Lemma 7 in the SM, 
we have $E_{j_1}\{M_{j_1j_2}^{2}\mathbb{I}(M_{j_1j_2} > \lambda_{p}(s))\} = E(L_{j_1j_2}^{2} | H_{0})\{1 + o_p(1)\}$, which implies
$E_{\mathscr{F}_{m-1}}\tilde{A}_{j_1j_2}^{2} = \text{Var}\{A_{j_1j_2} | H_{0}\}(1 + o_p(1))$, where $L_{j_1j_2} = M_{j_1j_2}\mathbb{I}(M_{j_1j_2} > \lambda_{p}(s))$.

Let $\rho_{j_1j_2 1} = \mbox{Cor}(X_{kj_1}, X_{kj_2})$ and $\rho_{j_1j_2 2} = \mbox{Cor}(Y_{kj_1}, Y_{kj_2})$ be the correlations. Let $\tilde{\rho}_{j_1j_2 1} = \tilde{\sigma}_{j_1j_2 1} / (\tilde{\sigma}_{j_1j_1 1}\tilde{\sigma}_{j_2j_2 1})^{1/2}$ and $\tilde{\rho}_{j_1j_2 2} = \tilde{\sigma}_{j_1j_2 2} / (\tilde{\sigma}_{j_1j_1 2}\tilde{\sigma}_{j_2j_2 2})^{1/2}$.
For $j_1 \neq j_3$ and $j_2 = j_4$,
by Lemma 7 in the SM,
$$\big| \text{Cor}_{(j_1,j_3)}(\tilde{\sigma}_{j_1j_2 1} - \tilde{\sigma}_{j_1j_2 2}, \tilde{\sigma}_{j_3j_2 1} - \tilde{\sigma}_{j_3j_2 2}) \big| \leq \tilde{\rho}_{j_1j_3},$$ where $\tilde{\rho}_{j_1j_3} = \max\{|\tilde{\rho}_{j_1j_3 1}|, |\tilde{\rho}_{j_1j_3 2}|\}$.
By Lemmas 6 and 7 in the SM, 
we have
\bea
|E_{(j_1,j_3)}(\tilde{A}_{j_1j_2}\tilde{A}_{j_3j_2})|
&\leq& L_p\tilde{\rho}_{j_1j_3} p^{-\frac{4s}{1 + \tilde{\rho}_{j_1j_3}}} \{1 + o_{p}(1)\} + O_{p}(L_p p^{-4s}n^{-1/2}). \nn
\eea
Similarly, for $j_1 = j_3$ and $j_2 \neq j_4$, by Lemma 7, we have that $\big| \text{Cor}_{j_1}(\tilde{\sigma}_{j_1j_2 1} - \tilde{\sigma}_{j_1j_2 2}, \tilde{\sigma}_{j_1j_4 1} - \tilde{\sigma}_{j_1j_4 2}) \big| \leq \rho_{j_2j_4}$ and $$|E_{j_1}(\tilde{A}_{j_1j_2}\tilde{A}_{j_1j_4})| \leq L_{p}\rho_{j_2j_4} p^{-{4s} / (1 + \rho_{j_2j_4})} \{1 + o_{p}(1)\} + O_{p}(L_{p} p^{-4s}n^{-1/2}),$$ where $\rho_{j_2j_4} = \max\{|\rho_{j_2j_4 1}|, |\rho_{j_2j_4 2}|\}$.
For $j_1 \neq j_3$ and $j_2 \neq j_4$, we have $$\big| \text{Cor}_{(j_1,j_3)}(\tilde{\sigma}_{j_1j_2 1} - \tilde{\sigma}_{j_1j_2 2}, \tilde{\sigma}_{j_3j_4 1} - \tilde{\sigma}_{j_3j_4 2}) \big| \leq \tilde{\rho}_{j_1j_3}\rho_{j_2j_4}.$$
By Assumption \ref{as5} and Davydov's inequality, for any positive constant $M$, there exists a constant $c > 0$ such that
\be
|E_{(j_1,j_3)}(\tilde{A}_{j_1j_2}\tilde{A}_{j_3j_4})| \leq C \gamma_{2}^{|j_2 - j_4|} \leq p^{-M}
\label{eq:LargeDist}\ee
for a constant $\gamma_2 \in (0, 1)$ and $|j_2 - j_4| > c \log p$.
For $j_2$ and $j_4$ close, by Lemmas 6 and 7, it follows that
$$|E_{(j_1,j_3)}(\tilde{A}_{j_1j_2}\tilde{A}_{j_3j_4})| \leq L_p \tilde{\rho}_{j_1j_3}\rho_{j_2j_4} p^{-{4s} / (1 + \tilde{\rho}_{j_1j_3}\rho_{j_2j_4})} \{1 + o_{p}(1)\} + O_{p}(L_{p} p^{-4s}n^{-1/2}).$$

Combining all the different cases above for the indexes $(j_1,j_2,j_3,j_4)$ together, equation (\ref{eq:MDiffVar}) can be decomposed as
\bea
E_{\mathscr{F}_{m-1}}(D_{m, 1}^{2}) &=& a^{2}(m-1)\text{Var}\{A_{12} | H_{0}\}\{1 + o_{p}(1)\} \nn \\
&+& \sum_{m_1,m_2 = 1}^{m-1} \sum_{j_1\in S_{m_1} \neq j_3 \in S_{m_2}} \sum_{j_2 \in S_{m}} E_{(j_1,j_3)}(\tilde{A}_{j_1j_2}\tilde{A}_{j_3j_2}) \label{eq:MDiffVar1} \\
&+& \sum_{m_1 = 1}^{m-1} \sum_{j_1\in S_{m_1}} \sum_{j_2 \neq j_4 \in S_{m}} E_{j_1}(\tilde{A}_{j_1j_2}\tilde{A}_{j_1j_4}) \label{eq:MDiffVar2} \\
&+& \sum_{m_1,m_2 = 1}^{m-1} \sum_{j_1\in S_{m_1} \neq j_3 \in S_{m_2}} \sum_{j_2 \neq j_4 \in S_{m}}
E_{\mathscr{F}_{m-1}}(\tilde{A}_{j_1j_2}\tilde{A}_{j_3j_4}). \label{eq:MDiffVar3}
\eea
Note that $\rho_{j_1j_3} = 0$ for $m_1 \neq m_2$ due to the independence between ${\bf Z}_{S_{m_1}}$ and ${\bf Z}_{S_{m_2}}$. 
Under Assumption \ref{as3}, we also have $|\tilde{\rho}_{j_1j_3} - \rho_{j_1j_3}| \leq L_p n^{-1/2}$ for $j_1\in S_{m_1}$ and $j_3\in S_{m_2}$.
The term in (\ref{eq:MDiffVar1}) is bounded by
\bea
a^3(m-1)^2 L_p n^{-1/2} p^{-4s}
+ a(m-1) \sum_{j_1\neq j_3 \in S_{m_1}} L_p \rho_{j_1j_3} p^{-\frac{4s}{1 + \rho_{j_1j_3}}}. \nn
\eea
By Assumption \ref{as5} and Davydov's inequality, for any $M > 0$, there exists a constant $c > 0$ such that $\rho_{j_1j_3} \leq p^{-M}$ for $|j_1 - j_3| > c \log p$.
Therefore, the summation of $L_p \rho_{j_1j_3} p^{-{4s} / (1 + \rho_{j_1j_3})}$ over $j_1 \neq j_3 \in S_{m_1}$ is bounded by
$$\sum_{|j_1 - j_3| \leq c \log p}L_p \rho_{j_1j_3} p^{-4s / (1 + \rho_{j_1j_3})} + \sum_{|j_1 - j_3| > c \log p} L_p p^{-M-4s} \leq aL_pp^{-4s / (2 - \epsilon)}$$
for a small positive constant $\epsilon > 0$.
For (\ref{eq:MDiffVar2}), similarly, we have
\vspace{-0.75cm}
\be\begin{split}
\bigg|\sum_{m_1 = 1}^{m-1} \sum_{j_1\in S_{m_1}} \sum_{j_2 \neq j_4 \in S_{m}} E_{j_1}(\tilde{A}_{j_1j_2}\tilde{A}_{j_1j_4})\bigg| &\leq a^3(m-1)O_{p}(L_p n^{-1/2} p^{-4s}) \\
+ \ a(m - 1) &\sum_{j_2 \neq j_4 \in S_{m}} L_p \rho_{j_2j_4} p^{-{4s} / (1 + \rho_{j_2j_4})},
\end{split}\nn\ee
which is bounded by $a^3(m-1)O_{p}(L_p n^{-1/2} p^{-4s}) + a^{2}(m - 1)L_p p^{-4s / (2 - \epsilon)}$.
For the last term in (\ref{eq:MDiffVar3}), by choosing $M$ in (\ref{eq:LargeDist}) sufficiently large, it is bounded by
$a^3(m-1)^2 L_p n^{-1/2} p^{-4s} + a^2(m-1) L_p p^{-4s / (2 - \epsilon)}$.

Notice that $\sigma^2_{m} = E_{\mathscr{F}_{m-1}}(D_{m, 1}^{2}) + O_{p}(L_p d^{2} a^{3} p^{-4s} n^{-1/2})$ by choosing $a \ll \sqrt{n}$.
Summing up all the terms in (\ref{eq:MDiffVar1}) -- (\ref{eq:MDiffVar3}), up to a multiplication of $1 + o_p(1)$, we have that
\bea
\sum_{m=1}^{d}\sigma^2_{m} &=& \frac{a^{2}d(d-1)}{2}\text{Var}\{A_{12} | H_{0}\} + O_{p}(p^3 L_p n^{-1/2} p^{-4s}) + O(p^2 L_pp^{-\frac{4s}{2 - \epsilon}}), \nn
\eea
where $a^{2}d(d-1)\text{Var}\{A_{12} | H_{0}\} / 2 
= \text{Var}(B_{1, n} | H_{0})(1 + o(1))$.
Since $L_pp^{2-\frac{4s}{2 - \epsilon}} = o \{\text{Var}(B_{1, n} | H_{0})\}$, it follows that
$$\sum_{m=1}^{d}\sigma^2_{m} = \text{Var}(B_{1, n} | H_{0})(1 + o(1)) + O_{p}(p^3 L_p n^{-1/2} p^{-4s}).$$
Note that $p^3 L_p n^{-1/2} p^{-4s} = o(L_pp^{2 - 2s})$ for any $n$ and $p$ when $s > 1/2$.
Given $n = p^{\xi}$ for $\xi\in(0, 2)$, $p^3 L_p n^{-1/2} p^{-4s}$ is at a small order of $\text{Var}(B_{1, n} | H_{0}) = L_pp^{2 - 2s}$ if $s > 1/2 - \xi/4$, which proves the first claim of (\ref{MCLT}).

For the second claim of (\ref{MCLT}), notice that $D_{m} = D_{m, 1} + D_{m, 2}$ where $|D_{m, 2}| \leq d L_p O(a^{2} p^{-2s} n^{-1/2})$.
Given $a \ll \sqrt{n}$, we have $\sum_{m=1}^{d}d^{4} a^{8} p^{-8s} n^{-2} \ll \text{Var}^{2}(B_{1, n} | H_{0})$ when $s > 1/4 - \xi/8$.
Since $D_{m}^{4} \leq 8 (D_{m, 1}^{4} + D_{m, 2}^{4})$, to show the second claim of (\ref{MCLT}),
we only need to focus on $D_{m, 1}^{4}$, which is
\bea
&& \sum_{m_1=1}^{m-1} \tilde{f}({\bf Z}_{S_{m_1}},{\bf Z}_{S_{m}})^4 + c_{2}\sum_{m_1, m_2}^{\ast} \tilde{f}({\bf Z}_{S_{m_1}},{\bf Z}_{S_{m}})^2 \tilde{f}({\bf Z}_{S_{m_2}},{\bf Z}_{S_{m}})^2 \nn \\
&& \ \ \ \ \ + \ c_{3}\sum_{m_1, m_2, m_3}^{\ast} \tilde{f}({\bf Z}_{S_{m_1}},{\bf Z}_{S_{m}})^2 \tilde{f}({\bf Z}_{S_{m_2}},{\bf Z}_{S_{m}}) \tilde{f}({\bf Z}_{S_{m_3}},{\bf Z}_{S_{m}}) \label{Dm4} \\
&& \ \ + \sum_{m_1, m_2, m_3, m_4}^{\ast} \tilde{f}({\bf Z}_{S_{m_1}},{\bf Z}_{S_{m}}) \tilde{f}({\bf Z}_{S_{m_2}},{\bf Z}_{S_{m}}) \tilde{f}({\bf Z}_{S_{m_3}},{\bf Z}_{S_{m}}) \tilde{f}({\bf Z}_{S_{m_4}},{\bf Z}_{S_{m}}) \nn
\eea
where $\sum^{\ast}$ indicates summation over distinct indices smaller than $m$, and $c_2$ and $c_3$ are two positive constants.

Note that the expectation of the last term in (\ref{Dm4}) equals to the expectation of its conditional expectation given ${\bf Z}_{S_{m}}$, where the conditional expectation is bounded by
\be\begin{split}
& \bigg| \sum_{m_1, m_2, m_3, m_4}^{\ast} \{E_{S_{m}}\tilde{f}({\bf Z}_{S_{m_1}},{\bf Z}_{S_{m}})\} \{E_{S_{m}}\tilde{f}({\bf Z}_{S_{m_2}},{\bf Z}_{S_{m}})\}\\
&~~~~ ~~~~~~~~~~~~~\times\{E_{S_{m}}\tilde{f}({\bf Z}_{S_{m_3}},{\bf Z}_{S_{m}})\} \{E_{S_{m}}\tilde{f}({\bf Z}_{S_{m_4}},{\bf Z}_{S_{m}})\} \bigg| \nn \\
\leq \ & \bigg( \sum_{m_1} | E_{S_{m}}\tilde{f}({\bf Z}_{S_{m_1}},{\bf Z}_{S_{m}}) | \bigg)^{4} = O\{a^{8}(m-1)^4 L_p p^{-8s} n^{-2}\}.
\end{split}\ee
Note that the summation of this quantity over $1 \leq m \leq d$ is a small order term of $\text{Var}^{2}(B_{1, n} | H_{0})$ giving $s > 1/4 - \xi/8$.

Next, following the derivation of $\sigma_{m}^{2}$, $E_{S_{m}} \tilde{f}({\bf Z}_{S_{m_1}},{\bf Z}_{S_{m}})^2$ is equal to
$$\sum_{j_1 \in S_{m_1}, j_2 \in S_{m}}\sum_{j_3 \in S_{m_1}, j_4 \in S_{m}} E_{S_{m}}(\tilde{A}_{j_1j_2}\tilde{A}_{j_3j_4}) = O\{L_pa^2 p^{-2s}\},$$
which leads to
\be\begin{split}
E_{S_{m}} \sum_{m_1, m_2}^{\ast} \tilde{f}({\bf Z}_{S_{m_1}},{\bf Z}_{S_{m}})^2 \tilde{f}({\bf Z}_{S_{m_2}},{\bf Z}_{S_{m}})^2
&= O\{a^4(m-1)^2L_pp^{-4s}\} \mbox{ \ and } \nn \\
E_{S_{m}} \sum_{m_1, m_2, m_3}^{\ast} \tilde{f}({\bf Z}_{S_{m_1}},{\bf Z}_{S_{m}})^2 \tilde{f}({\bf Z}_{S_{m_2}},{\bf Z}_{S_{m}}) &\tilde{f}({\bf Z}_{S_{m_3}},{\bf Z}_{S_{m}}) \nn \\
&= O\{a^{6} (m-1)^{3} L_p p^{-6s} n^{-1}\}. \nn
\end{split}\ee
The summation of the two terms above over $1 \leq m \leq d$ are at smaller orders of $\text{Var}^{2}(B_{1, n} | H_{0})$.
Also notice that
\bea
E \sum_{m_1=1}^{m-1} \tilde{f}({\bf Z}_{S_{m_1}},{\bf Z}_{S_{m}})^4
\leq \sum_{m_1=1}^{m-1} a^{6} \sum_{i \in S_{m_1} j \in S_{m}} E \tilde{A}_{ij}^{4} = a^{8}(m-1)L_pp^{-2s}.
\nn
\eea
Since $\sum_{m=1}^{d}a^{8}(m-1)L_pp^{-2s} = a^{6}L_pp^{2-2s} \ll p^{4-4s}$ if $a \ll p^{(1-s)/3}$, the second claim of (\ref{MCLT}) is valid given $a \ll \min\{n^{1/2}, p^{(1-s)/3}\}$ and $s > 1/4 - \xi/8$.
This proves the asymptotical normality of $T_{n}(s)$ for $s > 1/2 - \xi/4$ under $H_{0}$ of (\ref{H0}) by choosing $a \ll \min\{n^{1/2}, p^{(1-s)/3}\}$, $b \sim \max\{\log(p), \log(n_1+n_2)\}$ and $b \ll a$. $\square$


\begin{thebibliography}{11}

\newcommand{\enquote}[1]{``#1''}
\expandafter\ifx\csname
natexlab\endcsname\relax\def\natexlab#1{#1}\fi


\bibitem[{Anderson(2003)}]{Anderson_2003} Anderson, T. W. (2003).
\textit{An Introduction to Multivariate Statistical Analysis} (3rd
ed.), New York: John Wiley \& Sons.


\bibitem[{Arias-Castro et al.(2012)}]{AC_2012} Arias-Castro, E., Bubeck, S. and Lugosi, G. (2012).
Detection of correlations. \textit{The Annals of Statistics}, \textbf{40}, 412--435.

\bibitem[{Athreya and Lahiri(2006)}]{AthreyaLahiri} Athreya, K. and Lahiri, S. (2006).
\textit{Measure Theory and Probability Theory}, New York: Springer.

\bibitem[{Bai et al.(2009)}]{Bai_2009} Bai, Z. D., Jiang, D. D., Yao, J. F. and  Zheng S. R. (2009).
Corrections to LRT on large-dimensional covariance matrix by RMT. \textit{The Annals of Statistics}, \textbf{37}, 3822-3840.

\bibitem[{Bai and Silverstein(2010)}]{Bai_Silverstein_2010} Bai, Z. D. and Silverstein, J. W. (2010).
\textit{Spectral Analysis of Large Dimensional Random Matrices}, New York: Springer.


\bibitem[Bai and Yin(1993)]{Bai_Yin_1993} Bai, Z. D. and Yin, Y. Q. (1993). Limit of the smallest eigenvalue of a large dimensional sample covariance matrix.  \textit{The Annals of Probability}, \textbf{21}, 1275-1294.




\bibitem[{Berbee(1979)}]{Berbee_1979} Berbee, H. (1979).
\textit{Random Walks with Stationary Increments and Renewal Theory}, Amsterdam: Mathematical Centre.

\bibitem[Bickel and Levina(2008a)]{BL_2008a}
Bickel, P. and Levina, E. (2008a).
Regularized estimation of large covariance matrices. \textit{The Annals of Statistics},
\textbf{36}, 199-227.

\bibitem[Bickel and Levina(2008b)]{BL_2008b}
Bickel, P. and Levina, E. (2008b).
Covariance regularization by thresholding. \textit{The Annals of Statistics},
\textbf{36}, 2577-2604.


\bibitem[Bosq(1998)]{Bosq_1998}
Bosq, D. (1998). \textit{Nonparametric Statistics for Stochastic Processes} (2nd ed.),  New York: Springer.

\bibitem[Bradley(2005)]{Bradley_2005}
Bradley, R. (2005).
Basic properties of strong mixing conditions: a survey and some open questions.
\textit{Probability Surveys},
\textbf{2}, 107-144.



\bibitem[Cai et al.(2013)]{Cai_Liu_Xia_2013} Cai, T., Liu, W. D. and Xia, Y. (2013). Two-sample covariance matrix testing and support recovery in high-dimensional and sparse settings. \textit{Journal of the Americain Statistical Association}, \textbf{108}, 265-277.


\bibitem[Chang et al.(2017)]{Chang_2017} Chang, J. Y., Zhou, W.,  Zhou, W. X. and Wang, L. (2017). Comparing large covariance matrices under weak conditions on the dependence structure and its application to gene clustering. \textit{Biometrics}, \textbf{73}, 31-41.





\bibitem[{Delaigle et~al.(2011)}]{Delaigle_2011_JRSSB}
Delaigle, A., Hall, P. and Jin, J. (2011). Robustness and accuracy
of methods for high dimensional data analysis based on Student's
t-statistic. \textit{Journal of the Royal Statistical Society: Series
B (Statistical Methodology)}, \textbf{73}, 283-301.


\bibitem[{de la Fuente(2010)}]{Fuente_2010}
de la Fuente, A. (2010). From differential expression to differential networking--identification of dysfunctional regulatory networks in diseases. \textit{Trends in Genetics}, \textbf{26}, 326-333.


\bibitem[Donoho and Jin(2004)]{Donoho_Jin_2004}
Donoho, D. and Jin, J. (2004).
Higher criticism for detecting sparse heterogeneous mixtures.
\textit{The Annals of Statistics},
\textbf{32}, 962 - 994.

\bibitem[Donoho and Jin(2015)]{Donoho_Jin_2015}
Donoho, D. and Jin, J. (2015). Higher criticism for large-scale inference, especially for rare and weak effects.
\textit{Statistical Science}, \textbf{30}, 1--25.


\bibitem[Fan(1996)]{Fan_1996} Fan, J. (1996). Test of significance based on wavelet thresholding and Neyman's truncation.
\textit{Journal of the American Statistical Association}, \textbf{91}, 674-688.



\bibitem[Gupta and Giri(1973)]{Gupta_Giri_1973} Gupta, D. S. and Giri, N. (1973). Properties of tests concerning
covariance matrices of normal distributions. \textit{The Annals of
Statistics}, \textbf{6}, 1222-1224.

\bibitem[Hall and Heyde(1980)]{Hall_Heyde_1980}
Hall, P. and Heyde, C. C. (1980). \textit{Martingale Limit Theory and Its Application}, Academic Press.

\bibitem[Hall and Jin(2010)]{Hall_Jin_2010}
Hall, P. and Jin, J. (2010). Innovated higher criticism for
detecting sparse signals in correlated noise. \textit{The Annals of
Statistics}, \textbf{38}, 1686-1732.

\bibitem[Hotelling(1931)]{Hotelling_1931}
Hotelling, H. (1931). The generalization of student's ratio. \textit{Annals of Mathematical Statistics}, \textbf{2}, 54-65.


\bibitem[{Ingster(1997)}]{Ingster_1997}
Ingster, Y. I. (1997). Some problems of hypothesis testing leading
to infinitely divisible distributions. \textit{Mathematical Methods
of Statistics}, \textbf{6}, 47-69.



\bibitem[John(1971)]{John_1971}
John, S. (1971). Some optimal multivariate tests. \textit{Biometrika}, \textbf{59}, 123-127.


\bibitem[{Li and Chen(2012)}]{Li_Chen_2012} Li, J. and Chen, S. X. (2012). Two sample tests for
high-dimensional covariance matrices. \textit{The Annals of
Statistics}, \textbf{40}, 908-940.


\bibitem[Liu(2013)]{Liu_2013}
Liu, W. (2013). Gaussian graphical model estimation with false discovery rate control. \textit{The Annals of
Statistics}, {\bf 41}, 2948--2978.

\bibitem[Mokkadem(1988)]{Mokkadem_1988}
Mokkadem, A. (1988). Mixing properties of ARMA processes. \textit{Stochastic processes and their applications}, {\bf 29}, 309--315.

\bibitem[{Nagao(1973)}]{Nagao_1973} Nagao, H. (1973). On some test criteria for covariance matrix. \textit{The Annals of Statistics}, \textbf{1}, 700-709.


\bibitem[{Perlman(1980)}]{Perlman_1980} Perlman, M. D. (1980). Unbiasedness of the likelihood ratio tests for equality of several
covariance matrices and equality of several multivariate normal
populations. \textit{The Annals of Statistics}, \textbf{8}, 247-263.


\bibitem[{Qiu and Chen(2012)}]{QC_2012}
Qiu, Y. and Chen, S. X. (2012). Test for bandedness of high-dimensional covariance matrices and bandwidth estimation. \textit{The Annals of Statistics}, \textbf{40}, 1285-1314.

\bibitem[{Qiu et al.(2018)}]{Qiu_Chen_Nettleton_2016}
Qiu, Y., Chen, S. X. and Nettleton, D. (2018). Detecting rare and faint signals via thresholding maximum likelihood estimators. \textit{The Annals of Statistics},  \textbf{46}, 895-923.

\bibitem[Ren et al.(2015)]{Ren_2015}
Ren, Z., Sun, T., Zhang, C. H. and Zhou, H. (2015). Asymptotic normality and optimalities in estimation of large Gaussian graphical models. \textit{The Annals of Statistics}, {\bf 43}, 991--1026.

\bibitem[{Rothman(2012)}]{Rothman_2012} Rothman, A. J. (2012).
Positive definite estimators of large covariance matrices.
\textit{Biometrika}, \textbf{99}, 539-550.




\bibitem[{Saulis and Statulevi\v cius(1991)}]{Saulis_1991} Saulis, L. and Statulevi\v cius, V. A. (1991).  \textit{Limit Theorems for Large Deviations,} Dordrecht: Kluwer Academic.


\bibitem[{Schott(2007)}]{Schott_2007}
Schott, J. R. (2007). A test for the equality of covariance matrices
when the dimension is large relative to the sample sizes.
\textit{Computational Statistics and Data Analysis}, \textbf{51},
6535-6542.






\bibitem[{Srivastava and Yanagihara(2010)}]{Srivastava_Yanagihara_2010} Srivastava, M. S., and Yanagihara, H. (2010). Testing the equality
of several covariance matrices with fewer observations than the
dimension. \textit{Journal of Multivariate Analysis}, \textbf{101},
1319-1329.

%


\bibitem[Tran(1990)]{Tran_1990}
Tran, L. T. (1990). Kernel density estimation on random fields. \textit{Journal of Multivariate Analysis}, \textbf{34}, 37-53.



\bibitem[{Xue et al.(2012)}]{Xue_Ma_Zou_2012} Xue, L. Z., Ma, S. Q. and Zou, H. (2012). Positive-definite $\ell_1$-penalized estimation of large covariance matrices. \textit{Journal of the Americain Statistical Association}, \textbf{107},
1480-1491.

\bibitem[{Yi et al.(2007)}]{Yi_2007} Yi, G., Sze, S. H. and Thon, M. R. (2007). Identifying clusters of functionally related genes in genomes. \textit{Bioinformatics}, \textbf{23},
 1053-1060.	

\bibitem[{Zhong et al.(2013)}]{Zhong_Chen_Xu_2013} Zhong, P. S., Chen, S. X. and Xu, M. Y. (2013). Tests alternative to higher
criticism for high dimensional means under sparsity and column-wise
dependence. \textit{The Annals of Statistics}, \textbf{41},
2820-2851.



\end{thebibliography}
\end{document}